\documentclass[12pt]{amsart}

\usepackage{amssymb}

\usepackage{enumitem}

\usepackage{graphicx}

\makeatletter
\@namedef{subjclassname@2020}{%
	\textup{2020} Mathematics Subject Classification}
\makeatother

\usepackage[T1]{fontenc}

\newtheorem{theorem}{Theorem}[section]
\newtheorem{corollary}[theorem]{Corollary}
\newtheorem{lemma}[theorem]{Lemma}
\newtheorem{proposition}[theorem]{Proposition}

\theoremstyle{definition}
\newtheorem{definition}[theorem]{Definition}
\newtheorem{remark}[theorem]{Remark}

\numberwithin{equation}{section}

\DeclareGraphicsExtensions{.eps}
\DeclareFontFamily{U}{mathx}{}
\DeclareFontShape{U}{mathx}{m}{n}{<-> mathx10}{}
\DeclareSymbolFont{mathx}{U}{mathx}{m}{n}
\DeclareMathAccent{\widehat}{0}{mathx}{"70}
\DeclareMathAccent{\widecheck}{0}{mathx}{"71}
\allowdisplaybreaks

\begin{document}
	
	
	\baselineskip=17pt
	
	
	\title[Korovkin-type Approximation]{Korovkin-type Approximation for Non-positive Operators}
	
	\author[V. B. Kiran Kumar]{V. B. Kiran Kumar}
	\address{Cochin University of Science and Technology\\ Kalamassery\\
		Cochin\\
		India 682022}
	\email{vbk@cusat.ac.in}
	
	\author[M. N. N. Namboodiri]{M.N.N. Namboodiri}
	\address{IMRT Thiruvanathapuram and formerly of Cochin University of Science and Technology\\ Kalamassery\\ Cochin\\ India 682022}
	\email{mnnadri@gmail.com}
	
	\author[P. C. Vinaya]{P. C. Vinaya}
	\address{Cochin University of Science and Technology\\ Kalamassery\\
		Cochin\\
		India 682022}
	\email{vinayapc01@gmail.com}

	\date{}

	\begin{abstract}
	The classical Korovkin theorem traditionally relies on the positivity of the underlying sequence of operators. In 1968, D. E. Wulbert obtained a non-positive version by exploiting geometric properties of function spaces, namely the Choquet boundary and the unique extension property of extreme points of the dual unit ball for weakly separating subspaces. In this article we develop this geometric approach further and prove a Korovkin-type theorem for uniformly bounded sequences of operators on $C(X)$ and on $L^1[0,1]$, with the convergence of the operators on a test set being replaced by convergence to a limit operator. The main emphasis is on the underlying geometry: the weakly separating subspace, its Choquet boundary, and the associated unique extension property. As an application, we show that the Grünwald interpolation operators, which are non-positive, satisfy a Korovkin-type approximation theorem. We also extend this to an $L^1(\mathbb{R})$ setting and include numerical illustrations.
	\end{abstract}
	
	
	\keywords{Korovkin theorem, Wulbert's theorem, Operator Version of Korovkin Theorem, Non-positive Operators, Modulus of continuity}
	
	\maketitle
	
\section{Introduction}
In 1953, P. P. Korovkin established the following fundamental result in approximation theory, now known as Korovkin's theorem. Let $\| \cdot \|_\infty$ denote the sup-norm on $[a,b]$. 
\begin{theorem}\cite{korovkin}
	Let $\{L_n\}_{n\in\mathbb{N}}$ be a sequence of positive linear operators on $C[a,b]$ such that 
	$\lim\limits_{n\to\infty}\|L_n(g)-g\|_\infty=0$ for $g\in\{1,x,x^2\}$. Then 
	$\lim\limits_{n\to\infty}\|L_n(f)-f\|_\infty=0$ for every $f\in C[a,b]$. 
\end{theorem}
In other words, convergence on the three test functions $1$, $x$, $x^2$ already forces convergence on the whole space $C[a,b]$. This result unified several classical approximation processes such as those of Bernstein and Weierstrass. The periodic analogue on $C_{2\pi}(\mathbb{R})$ uses the test functions $1$, $\sin x$, $\cos x$. Numerous extensions of Korovkin theorem have been obtained, most of them relying on positivity (or complete positivity in the non-commutative setting) of the operators; see, for instance, \cite{altomarecampiti}.

A decisive step away from positivity was taken by D. E. Wulbert in 1968. He proved the following two results.

\begin{theorem}[Corollary $3$, \cite{wulbert}]\label{wul1}
	Let $\{L_n\}_{n\in\mathbb{N}}$ be a sequence of norm-one operators defined on $C[0, 1]$
	($C_{2\pi}(\mathbb{R})$, respectively). Then $L_n(f)$ converges to $f$ for all $f
	\in C [0,1]$ ($C_{2\pi}(\mathbb{R})$), if and
	only if $L_n(p)$ converges to $p$ for the three functions $1$, $x$ and 
	$x^2$ ($1$, $\cos x$ and $\sin x$, respectively).
\end{theorem}

This result has been extended to the space $L^1[0,1]$, as stated below. We denote $\|L\|$ as the operator norm of an operator $L$.
\begin{theorem}[Theorem $9$ of \cite{wulbert}]\label{wul2}
	Let $\{L_n\}_{n\in\mathbb{N}}$ be a sequence of operators defined on $L^1[0, 1]$. If
	\begin{enumerate}
		\item $L_n(1)$ converges to $1$,
		\item $L_n(p)$ converges weakly to $p$ for the two functions $p = x$ and $p= x^2$,
		\item $\lim\limits_{n \to \infty}  \|L_{n}\| = 1$,
	\end{enumerate}
	then $L_n(f)$ converges to $f$ for all $f$ in $L^1[0, 1]$.
\end{theorem}

Wulbert's theorems are remarkable because they do not assume positivity. Instead, the proofs rely on a geometric mechanism. The key ingredients are the concepts of a weakly separating subspace and the associated Choquet boundary. Let $P$ be a linear subspace of a normed space $E$. We say that $P$ is weakly separating if every extreme point of the unit ball of $P^*$ has a unique extension to an extreme point of the unit ball of $E^*$. The set of all such extended extreme points is the Choquet boundary of $P$, denoted by $\operatorname{ch}(P)$. This unique extension property is a purely geometric condition and replaces the usual order-theoretic positivity condition.

The geometric approach to Korovkin-type theorems has been further developed in various articles, in commutative as well as non-commutative settings \cite{bvlshirali, mnnrims, mnnoam}. These works highlight the role of the Choquet boundary and the unique extension property in obtaining convergence theorems for sequences of operators.

The purpose of the present article is to continue this geometric line of investigation. Our main results, Theorems \ref{thm1} and \ref{thm2}, generalize Wulbert's theorems in two directions. First, we replace the requirement that the operators have norm one (or that their norms converge to one) by the weaker assumption of uniform boundedness. Second, we allow convergence to a bounded linear operator $L$ rather than only to the identity operator. The proofs remain essentially geometric. The crucial step is Lemma \ref{lemm1}, which is a modification of Wulbert's lemma that we need in order to establish the convergence of the Gr\"unwald interpolation operators via the Korovkin-type theorem.

As a concrete illustration, we consider the Gr\"unwald interpolation operators, $\{G_n\}_{n\in\mathbb{N}}$ on $C[0,\pi]$. These operators are uniformly bounded but not positive. Using the geometric Korovkin-type theorem, we show that $G_n(f)$ converges to $f$ for every $f\in C[0,\pi]$. We also construct an extension of the Gr\"unwald operators to $L^1(\mathbb{R})$ via the Fourier transform and an approximate identity, obtaining a family $\{H_{n,\delta}\}_{n\in\mathbb{N}, \delta>0}$ that converges on compact intervals. In the final part of the article we present a direct approach to the convergence of $\{G_n\}_{n\in\mathbb{N}}$ and $\{H_{n,\delta}\}_{n\in\mathbb{N},\delta>0}$, including numerical illustrations and a rate of convergence estimate expressed in terms of the modulus of continuity.

The article is organized as follows. In Section \ref{wulb} we prove a slightly general version of Wulbert's theorem on $C(X)$ and on $L^1[0,1]$.We also recall the Gr\"unwald operators, their non-positivity, and extend the notion to $L^1(\mathbb{R})$, using the classical Fourier transform techniques. Section \ref{direct} contains the direct convergence proof, the rate of convergence estimate, and numerical tables. We conclude with some remarks and open problems.
\section{An Operator Version of Wulbert's Theorem}\label{wulb}

In this section, we prove an operator version of Theorem \ref{wul1} and \ref{wul2}.\par
First, we recall some preliminary definitions and notations from \cite{wulbert}.\par

Let $E$ be a normed linear space, and $E^*$ be the dual space of $E.$ Denote  $\mathbb{S}_{E^*}$ to the closed unit ball of the space $E^*$ and $\mathrm{ext }\ \mathbb{S}_{E^*}$, the set of all extreme points of $\mathbb{S}_{E^*}$. The weak topology on $E^*$ induced by $E$ is denoted by $w(E^*,E)$.
We recall the definition of a weakly separating subspace and the generalized Choquet boundary given in \cite{wulbert}.
\begin{definition}
	Let $P$ be a subspace of a normed linear space $E$. Then $P$ is a weakly separating subspace of $E$ if for each element $L\in \mathrm{ext}\ \mathbb{S}_{P^*}$, there is a unique element $\Phi$ of $\mathrm{ext}\ \mathbb{S}_{E^*}$ such that $\Phi|_P=L$.
\end{definition}
\begin{definition}
	Let $P$ be weakly separating in $E$. The generalized Choquet boundary of $P$ is defined to be the set;
	\[
	\mathrm{ch(P)}=\{F\in \mathrm{ext}\ \mathbb{S}_{E^*}: F|_P\in \mathrm{ext}\ \mathbb{S}_{P^*}\}.
	\]
\end{definition}
We denote the closed ball of radius $\lambda>0$ in $E^*$ by $\mathbb{S}^\lambda_{E^*}$, which is compact with respect to the topology $w(E^*, E)$ (Banach Alaoglu theorem). We also note that $\phi\in \mathbb{S}^\lambda_{E^*}$ is an extreme point of $\mathbb{S}^\lambda_{E^*}$ if and only if $\frac{\phi}{\lambda}$ is an extreme point of $\mathbb{S}_{E^*}$.\par 
We extend the notion of generalized Choquet boundary as follows. The following lemma ensures that the extension is well-defined.
\begin{lemma}\label{lemm0}
	Suppose $P$ is weakly separating in $E$. Then for each element $L'\in \mathrm{ext}\ \mathbb{S}^\lambda_{P^*}$, there exists a unique member $\Phi\in \mathrm{ext}\ \mathbb{S}^\lambda_{E^*}$ such that $\Phi|_P=L'$.
\end{lemma}
\begin{proof}
	Let $L'\in \mathrm{ext}\ \mathbb{S}^\lambda_{P^*}$. Then $\frac{L'}{\lambda}\in \mathrm{ext}\ \mathbb{S}_{P^*}$. Since $P$ is weakly separating in $E$, there exists a unique $F\in \mathrm{ext}\ \mathbb{S}_{E^*}$ such that $F|_P=\frac{L'}{\lambda}$. Define $\phi=\lambda F$. Then $\phi\in \mathrm{ext}\ \mathbb{S}^\lambda_{E^*}$ and $\phi|_P=L'$. The uniqueness of $\phi$ follows from the uniqueness of $F$.
\end{proof}

\begin{definition}\label{def:chlambda}	
	Suppose $P$ is weakly separating in $E$. The generalized Choquet Boundary of $P$ with respect to $\mathbb{S}^\lambda_{E^*}$ is defined by,
	\[
	\mathrm{ch}_\lambda(P)=\{\phi \in \mathrm{ext} \ \mathbb{S}^\lambda_{E^*}: \phi|_P\in \mathrm{ext}\  \mathbb{S}^\lambda_{P^*} \}.
	\]
\end{definition}
We note that $\mathrm{ch}_\lambda(P)=\lambda \mathrm{ch}(P)$.

\begin{remark}\label{rem:chlambda}
	The scaled Choquet boundary $\mathrm{ch}_\lambda(P)$ provides a natural framework for stating the main lemma. However, for the proof, it is more convenient to work with the unscaled version $\mathrm{ch}(P)$ and handle the scaling factor $\lambda$ explicitly. Throughout this article, we shall use $\mathrm{ch}(P)$ in the statements and proofs of the main results, while $\mathrm{ch}_\lambda(P)$ serves as a conceptual bridge.
\end{remark}

The following lemma is a generalization of Lemma $1$ in \cite{wulbert}. This lemma plays a crucial role in obtaining our main results.
\begin{lemma}\label{lemm1}
	Let $P$ be a weakly separating subspace of a normed linear space $E$ and $M$ be a linear subspace of $E$ which contains $P$. Let $D$ be a dense subspace of $M$. Let $\{L_n\}$ be a net of norm one operators that carry $M$ into $E$ satisfying 
	\begin{equation}\label{eqbd1}
		\limsup\limits_{n}|\Phi(L_n(f))|\leq \lambda\|f\| 
	\end{equation}

	for all $f\in D$ and $\Phi\in \mathrm{ch}(P)$, where $0<\lambda\leq 1$. Suppose $f\in M$. 
	\begin{enumerate}
		\item If 
		$L_n(p)$ converges to $\lambda p$ in 
		$w(M^*,M)$ for all $p\in P$, then $\Phi\circ L_n(f)$ converges to $\lambda \Phi(f)$ for all $\Phi$ in $\mathrm{ch}(P)$.
		\item Suppose that the bound in \eqref{eqbd1} holds uniformly on all compact subsets of $\mathrm{ch}(P)$. If $L_n(p)$ converges to $\lambda p$ in norm for all $p$ in $P$, then $L_n(f)$ converges uniformly to $\lambda f$ on $w(E^*,E)$-compact subsets of $\mathrm{ch}(P)$.
	\end{enumerate}
\end{lemma}

\begin{proof}
	We prove $(1)$. Suppose we have 
	$\phi\circ L_n(p)$ converges to $\lambda \phi(p)$ for all $\phi \in M^*$ and for all $p\in P$. Let $\Phi\in \mathrm{ch}(P)$ and let $\{L_i\}$ be any subnet of the net $\{L_n\}$. \par
	We prove that this subnet $\{L_i\}$ has further a subnet say $\{L_j\}$ such that $\Phi\circ L_j$ converges in the $w(M^*,M)$-topology. Then we prove that the limit functional is $\lambda \Phi$. Hence we obtain $\Phi\circ L_n(f)$ converges to $\lambda \Phi(f)$ for all $\Phi\in \mathrm{ch}(P)$ and we are done.\par 
	By hypothesis $\|L_i\|=1$. Hence we have $\Phi\circ L_i \in \mathbb{S}_{M^*}$. But we know that $\mathbb{S}_{M^*}$ is compact in the $w(M^*,M)$-topology by the well known Banach Alaoglu Theorem. Thus $\Phi\circ L_i$ has further a subnet $\Phi\circ L_j$ which converges to some functional say $H$. Also by the hypothesis, for all $g\in D$, we have 
	\begin{align*}
		|H(g)|&=\lim\limits_{j}|\Phi\circ L_j(g)|\\
		&=\limsup\limits_{j}|\Phi\circ L_j(g)|\\
		&\leq \lambda \|g\|.
	\end{align*}
	By the density of $D$, we also have $|H(f)|\leq \lambda\|f\|$ for all $f\in M$. Therefore, $\|H\|\leq\lambda$, so $\frac{H}{\lambda}\in \mathbb{S}_{M^*}\subseteq\mathbb{S}_{E^*}$ (by Hahn-Banach extension).
	We also have 
	\[
	H(p)=\lambda \Phi(p)\ \text{for all}\ p\in P,
	\]
	so that $\frac{H}{\lambda}|_P = \Phi|_P$. Since $\Phi\in \mathrm{ch}(P)$, we have $\Phi|_P\in \mathrm{ext}\,S_{P^*}$.
	
	We now establish the compound fact that $\frac{H}{\lambda} = \Phi$ on $E$. Define the face
	\[
	C := \{F\in\mathbb{S}_{E^*} : F|_P = \Phi|_P\}.
	\]
	This is a non-empty, closed (in $w(E^*,E)$), convex subset of $\mathbb{S}_{E^*}$, hence compact by Banach-Alaoglu.
	
	\textbf{Claim 1:} $\mathrm{ext}\,C \subseteq \mathrm{ext}\,\mathbb{S}_{E^*}$.
	
	\textit{Proof of Claim 1:} Let $F\in\mathrm{ext}\,C$. Suppose $F = \frac{F_1+F_2}{2}$ with $F_1,F_2\in\mathbb{S}_{E^*}$. Restricting to $P$, we get:
	\[
	\Phi|_P = F|_P = \frac{F_1|_P+F_2|_P}{2}.
	\]
	Since $\Phi|_P\in\mathrm{ext}\,\mathbb{S}_{P^*}$, we must have $F_1|_P = F_2|_P = \Phi|_P$. Thus $F_1,F_2\in C$. Since $F\in\mathrm{ext}\,C$, we get $F_1=F_2=F$. This proves $F\in\mathrm{ext}\,\mathbb{S}_{E^*}$.
	
	\textbf{Claim 2:} $\mathrm{ext}\,C = \{\Phi\}$.
	
	\textit{Proof of Claim 2:} By Claim 1, $\mathrm{ext}\,C \subseteq \mathrm{ext}\,\mathbb{S}_{E^*}$. By the weak separation property, $\Phi$ is the unique element of $\mathrm{ext}\,\mathbb{S}_{E^*}$ that agrees with $\Phi|_P$ on $P$. Therefore, $\mathrm{ext}\,C = \{\Phi\}$.
	
	\textbf{Claim 3:} $C = \{\Phi\}$.
	
	\textit{Proof of Claim 3:} By the Krein-Milman theorem, $C = \overline{\mathrm{conv}}(\mathrm{ext}\,C) = \overline{\mathrm{conv}}(\{\Phi\}) = \{\Phi\}$.
	
	Since $\frac{H}{\lambda}\in C$, we conclude $\frac{H}{\lambda} = \Phi$ on $E$, hence $H = \lambda\Phi$ on $M$. This proves $(1)$.
	\par
	\par 
	We prove $(2)$ by the method of contradiction as in \cite{wulbert}. Let $K$ be a compact subset of $\mathrm{ch}(P)$ and suppose $L_n(f)$ does not converge uniformly to $\lambda f$ on $K$. Then there exists a subnet $\{L_j\}$ of $\{L_n\}$ and a sequence $\Phi_j$ in $K$ and some $\epsilon>0$ such that 
	\[
	|\Phi_j\circ L_j(f)-\lambda \Phi_j(f)|\geq \epsilon.
	\]
	We have that $\Phi_j\in \mathrm{ch}(P)$. Thus each $\Phi_j\circ L_j$ has an extension to $\mathbb{S}_{E^*}$. Let $H_j$ be the extension. Since both $K$ and $\mathbb{S}_{E^*}$ are compact in the $w(E^*,E)$ topology, we can find another subnet of $\Phi_j$ and $H_j$ say; $\Phi_l$ and $H_l$ respectively such that this subnet converges to some $\Phi$ and $H$ respectively.
	
	We first establish the norm bound on $H$. For each $g\in D$, we have:
	\[
	|H(g)| = \lim_l |H_l(g)| = \lim_l |\Phi_l\circ L_l(g)| \leq \limsup_l |\Phi_l\circ L_l(g)| \leq \lambda\|g\|,
	\]
	where in the last inequality we use the hypothesis uniformly over $\Phi_l\in K\subset\mathrm{ch}(P)$. By density of $D$ in $M$, we conclude $\|H\|\leq\lambda$, so $\frac{H}{\lambda}\in\mathbb{S}_{E^*}$.
	
	Let $p\in P$. Consider
	\[
	|\lambda \Phi(p)- H(p)|\leq |\lambda \Phi(p)-\lambda \Phi_l(p)|+|\lambda \Phi_l(p)- H_l(p)|+ |H_l(p)-H(p)|.
	\]
	In the right hand side, the first and the third term converges to $0$ as already seen. For the second term we have 
	\[
	\lambda \Phi_l(p)- H_l(p)=\lambda \Phi_l(p)-\Phi_l\circ L_l(p)=\Phi_l(\lambda p-L_l(p)).
	\]
	By hypothesis, we have $\lim\limits_{l}\|L_l(p)-\lambda p\|=0,$
	where $\|.\|$ denotes the norm on $E$. Thus the second term also converges to $0$. Hence we have $H(p)=\lambda \Phi(p)$ for all $p\in P$, so $\frac{H}{\lambda}|_P = \Phi|_P$.
	
	Since $\Phi\in K\subset\mathrm{ch}(P)$ (as $K$ is compact, hence closed), we have $\Phi|_P\in\mathrm{ext}\,\mathbb{S}_{P^*}$. By the same compound fact (Claims 1-3) applied to $C := \{F\in\mathbb{S}_{E^*} : F|_P = \Phi|_P\}$, we conclude $\frac{H}{\lambda} = \Phi$ on $E$, i.e., $H=\lambda\Phi$ on $E$.
	
	Also $H_l(f)$ converges to $H(f)=\lambda \Phi(f)$ for every $f\in E$. Now,
	\[
	|\Phi_l\circ L_l(f)-\lambda \Phi_l(f)|\leq |\Phi_l\circ L_l(f)-H(f)|+|\lambda \Phi(f)-\lambda \Phi_l(f)|
	\]
	The first term converges to $0$ since $\Phi_l\circ L_l(f)=H_l(f)$ converges to $H(f)=\lambda\Phi(f)$. The second term converges to $0$ since $\Phi_l$ converges to  $\Phi$ in the $w(E^*,E)$ topology. This contradicts our assumption. This proves $(2)$.
\end{proof}
\begin{remark}\label{eval}
	Let $X$ be a compact Hausdorff space. Let $P$ be
	a linear subspace of $C(X)$ which contains the constants and separates the
	points of $X$. Let $E=C(X)$, then $\mathrm{ch}(P)$ consists of all evaluation functionals and negative of evaluation functionals on the classical Choquet boundary of $X$ (see page $29$, \cite{phelps} and \cite{wulbert} for more details). 
\end{remark}
Now, we prove our main result.
\begin{theorem}\label{thm1}
	Let $P$ be a linear subspace of $C(X)$ containing constants such that $P$ separates points of $X$ and the classical Choquet boundary of $P$ is $X$. Let $M$ be a linear subspace of $C(X)$ which contains $P$. Suppose $\{L_n\}_{n\in\mathbb{N}}$ be a sequence of linear operators from $M$ onto $C(X)$ satisfying
	\begin{enumerate}
		\item $\sup\limits_{n}\|L_n\|=\alpha<\infty$\\
		\item $\limsup\limits_{n\to\infty}\|L_n(f)-L(f)+f\|_\infty\leq \|f\|_\infty$ for all $f\in D$, where $D$ is a dense subset of $M$.
	\end{enumerate}
	and $L$ be a bounded linear operator from $M$ onto $C(X)$. If $L_n(p)$ converges to $L(p)$ for all $p\in P$, then $L_n(f)$ converges to $L(f)$ for all $f\in M$. 
\end{theorem}
\begin{remark}
	The condition on the dense subset $D$ in Theorem \ref{thm1} is a technical hypothesis, inherited from Lemma \ref{lemm1}, which is used to control the norm of certain limit functionals arising in the proof. It is a standard assumption in this type of generalization.
\end{remark}
\begin{proof}
	Suppose that the limit operator is the identity operator, that is $L=I$. Note that $P$ contains the constant function $1$. We see that, whenever $\lim\limits_{n\to\infty}L_n(1)=1$,  $\lim\limits_{n\to\infty}\|L_n(1)\|_\infty=1$, $\sup\limits_{n}\|L_n\|\geq \sup\limits_{n}\|L_n(1)\|_\infty\geq 1$. Thus, we must have, $\sup\limits_{n}\|L_n\|=\alpha\geq 1$. \\
	\textbf{Step 1:}
	We claim the following:
	If $\lim\limits_{n\to\infty}\|L_n\|=\beta\geq 1$ and if $\lim\limits_{n\to \infty}\|L_n(p)-p\|_\infty= 0$ for all $p\in P$ then, $\lim\limits_{n\to \infty}\|L_n(f)-f\|_\infty= 0$ for all $f\in M$. \par 
	To prove this let $\lim\limits_{n\to \infty}\|L_n(p)-p\|_\infty=0$ for all $p\in P$. Suppose $\tilde{L}_n=\frac{L_{n}}{\|L_n\|}$ and let $\lambda=\frac{1}{\beta}\leq 1$. Then we have, $\lim\limits_{n\to \infty}\|\tilde{L}_n(p)-\lambda p\|_\infty= 0$ for all $p\in P$. 
	
	We now verify the hypotheses of Lemma \ref{lemm1} for the sequence $\{\tilde{L}_n\}$. For $f\in D$ and $\Phi\in\mathrm{ch}(P)$, using the condition on $D$ for the original sequence $\{L_n\}$ (with $L=I$), we have:
	\[
	\limsup_{n\to\infty}|\Phi(\tilde{L}_n(f))| = \limsup_{n\to\infty}\frac{|\Phi(L_n(f))|}{\|L_n\|} \leq \frac{1}{\beta}\|f\|_\infty = \lambda\|f\|_\infty,
	\]
	since $\lim\limits_{n\to\infty}\|L_n\|=\beta$ and $\limsup\limits_{n\to\infty}|\Phi(L_n(f))|\leq \|f\|_\infty$ (from the condition on $D$ with $L=I$). This verifies the hypothesis of Lemma \ref{lemm1}.
	
	Then part $(2)$ of Lemma \ref{lemm1} gives that $\tilde L_n(f)$ converges to $\lambda f$ on $w(C(X)^*,C(X))$ - compact subsets of $\mathrm{ch}(P)$.\par 
	Let $x\in X$, then the evaluation functional on $C(X)$ is the map $ev_x:C(X)\rightarrow \mathbb{C}$ defined by $ev_x(f)=f(x)$. Thus, we have $ev_x, -ev_x\in \mathrm{ch}(P)$ for all $x\in X$ (see Remark \ref{eval}). Since $X$ is compact, the set $\{ev_x:x\in X\}$ and $\{-ev_x:x\in X\}$ are compact in $\mathrm{ch}(P)$ with respect to the $w(C(X)^*,C(X))$ topology. Therefore by Lemma \ref{lemm1},
	\[
	|ev_x(\tilde L_n(f))-ev_x(\lambda f)|=|\tilde L_n(f)(x)-\lambda f(x)|,
	\]
	which converges to $0$ uniformly for all $x\in X$,
	which gives $\lim\limits_{n\to \infty}\|\tilde L_n(f)-\lambda f\|_\infty=0$. Hence, $\lim\limits_{n\to \infty}\|L_n(f)-f\|_\infty= 0$ for all $f\in M$. \\
	\textbf{Step 2:} If  $\sup\limits_{n}\|L_{n}\|=\alpha<\infty$, then $\lim\limits_{n\to \infty}\|L_n(p)-p\|_\infty=0$ for all $p\in P$ gives $\lim\limits_{n\to \infty}\|L_n(f)-f\|_\infty=0$ for all $f\in M$.\par 
	Suppose that $\lim\limits_{n\to \infty}\|L_n(p)-p\|_\infty=0$ for all $p\in P$.
	We need to prove that $\lim\limits_{n\to \infty}\|L_n(f)-f\|_\infty= 0$ for all $f\in M$.
	Suppose there exist an $f\in M$ such that this is not true. Then there exist a sub-sequence $\{L_{n_l}\}_{l\in\mathbb{N}}$ of $\{L_n\}_{n\in\mathbb{N}}$ and an $\epsilon>0$ such that 
	\begin{equation}\label{con}
		\|L_{n_l}(f)-f\|_\infty>\epsilon 
	\end{equation}
	for all $l\in\mathbb{N}$.
	Since $\sup\limits_n\|L_n\|=\alpha\geq 1$, we have $\sup\limits_{n_l}\|L_{n_l}\|=\alpha_1$, for some $\alpha_1\leq \alpha$. Note that $\alpha_1\geq1$. This follows since $\lim\limits_{l\to\infty}\|L_{n_l}(1)\|_\infty=1$ and $\sup\limits_{n_l}\|L_{n_l}\|\geq \sup\limits_{n_l}\|L_{n_l}(1)\|_\infty\geq 1$. Hence, we can further find a sub-sequence of $\{L_{n_l}\}_{l\in\mathbb{N}}$, say $\{L_{n_k}\}_{k\in\mathbb{N}}$ such that $\lim\limits_{k\to\infty}\|L_{n_k}\|=\alpha_1$, which gives $\lim\limits_{k\to\infty}\|L_{n_k}(f)-f\|_\infty=0$ for all $f\in M$. But this is a contradiction to \ref{con}, Therefore, we must have $\lim\limits_{n\to\infty}\|L_n(f)-f\|_\infty=0$ for all $f\in M$.\\
	\textbf{Step 3:} Let $\sup\limits_{n\to\infty}\|L_{n}\|<\infty$ and 
	suppose that $\lim\limits_{n\to\infty}\|L_n(p)-L(p)\|_\infty= 0$ for all $p\in P$. Let $\mathcal{L}_n=L_n+I-L$. Then, we have $\sup\limits_{n\to\infty}\|\mathcal{L}_n\|<\infty$ and we have $\lim\limits_{n\to\infty}\|\mathcal{L}_n(p)-p\|_\infty= 0$ for all $p\in P$. Hence by step $2$,  $\lim\limits_{n\to\infty}\|\mathcal{L}_n(f)-f\|_\infty=0$ for all $f\in M$.
\end{proof}
Now we apply Theorem \ref{thm1} to obtain results on spaces $C[0,1]$ and $C_{2\pi}(\mathbb{R})$. First, we recall the definition of a peak point.
\begin{definition}
	Let $X$ be a compact Hausdorff space and $P$ a closed subspace of $C(X)$, separating points and containing the constants. A point $x_0\in X$ is a peak point of $P$ if there exists a $g\in P$ for which $|g(x)|< |g(x_0)|$, for $x\neq x_0$.
\end{definition}
The following result connects the peak points and the Choquet boundary of a set.
\begin{theorem}\cite{phelps}\label{peak}
	Let $X$ be a compact Hausdorff space and let $P$ be the linear subspace spanned by a subset $S$ of $C(X)$ that contains constants and separates points of $X$. Then the set of peak points of $P$ is contained in the Choquet boundary of $P$.
\end{theorem}
\begin{corollary}\label{cor5.10}
	Let $\{L_n\}_{n\in\mathbb{N}}$ be a sequence of operators defined on $C[0, 1]$ ($C_{2\pi}(\mathbb{R})$, respectively) satisfying
	\begin{enumerate}
		\item $\sup\limits_{n}\|L_n\|=\alpha<\infty$\\
		\item $\limsup\limits_{n\to\infty}\|L_n(f)-L(f)+f\|_\infty\leq \|f\|_\infty$ for all $f\in D$, where $D$ is a dense subset of $C[0,1]$.
	\end{enumerate}
	and $L$ be a bounded linear operator on $C[0,1]$. Then $L_n(f)$ converges to $L(f)$ uniformly for all $f$ in $C[0, 1]$ ($C_{2\pi}(\mathbb{R})$, respectively) if and only if $L_n(p)$ converges to $L(p)$ for the three functions $1$, $x$, and $x^2$ ($1$, $\cos x$ and $\sin x$ respectively). 
\end{corollary}
\begin{proof}
	The Choquet boundary spanned by the set $\{1, x, x^2\}$ is $[0,1]$. Let $P=span\{1, x,x^2\}$ and $x_0\in [0,1]$. To see this,
	let $g(x)=1-(x-x_0)^2$, then $g\in P$. Now $|g(x)|<|g(x_0)|$ for all $x\in [0,1]$. Therefore, $x_0$ is a peak point of $P$. Thus all points in $[0,1]$ are peak points. Since peak points are contained in the Choquet boundary of $P$, the Choquet boundary spanned by $P$ is $[0,1]$. \par 
	Similarly let $P=span \{1, \sin x, \cos x\}$ and $x_0\in [0,2\pi]$. 
	Let $g(x_0)=1+\cos (x-x_0)=2-2\sin^2\frac{(x-x_0)}{2}$, then $g\in P$ and $|g(x)|<|g(x_0)|$ for all $x\in [0,2\pi]$ and hence the Choquet boundary spanned by the set $\{1, \sin x,\cos x\}$ is $[0, 2\pi]$. Now, the result follows by applying Theorem \ref{thm1}.
\end{proof}
Now we extend these results to the space $L^1[0,1]$. We recall the following result (Corollary $7$ from \cite{wulbert}). 
\begin{corollary}\cite{wulbert}\label{sgn}
	If $P$ is a finite dimensional subspace of $L^1[0,1]$ spanned by polynomials or by trigonometric polynomials, then $P$ is weakly separating in $L^1[0,1]$ and,
	\[
	\mathrm{ch}(P)=\{sgn(f):f\in P\}
	\]
	where $sgn$ denotes the signum function defined by $sgn(f)(x)=\frac{f(x)}{|f(x)|}$ for $f\in L^1[0,1]$, $f\neq 0$.
	
\end{corollary}
We denote by $\|\cdot\|_1$, the $L^1$ norm on an interval $[a,b]\subset\mathbb{R}$, $\|f\|_1=\int\limits_a^b |f(x)|\,dx .$ and $\chi_{[a,b]}$ denotes the characteristic function defined on $[a,b]$.
\begin{lemma}\label{lemm2}
	Let $\{L_n\}_{n\in\mathbb{N}}$ be a sequence of norm one operators defined on $L^1[0, 1]$ and let $0<\lambda\leq 1$. If
	\begin{enumerate}
		\item For each characteristic function $g=\chi_{[a,b]}$, where $[a,b]\subset [0,1]$,
		\[
		\limsup\limits_{n\to\infty}\|L_n(g)\|_1\leq \lambda\|g\|_1,
		\]
		\item $L_n(1)$ converges to $\lambda.1$,
		\item $L_n(p)$ converges weakly to $\lambda p$ for the two functions $p = x$ and $p= x^2$.
	\end{enumerate}
	Then $L_n(f)$ converges to $\lambda f$ for all $f$ in $L^1[0, 1]$.
\end{lemma}

\begin{proof}
	We begin as in the proof of Theorem $9$ in \cite{wulbert}. 
	
	\textbf{Step 1}: Let $P=span\{1,x,x^2\}$. Thus, $L_n(p)$ converges to  $\lambda p$ weakly for all $p\in P$. Let
	\[
	K=\{ g\in L^\infty[0,1]:g=\chi_{[a,b]}, \  [a,b]\subset [0,1]\}.
	\]
	By Corollary \ref{sgn}, $\mathrm{ch}(P)=\{sgn f: f\in P\}$. It is already observed in \cite{wulbert} that every member of $K$ can be written as the average of two members of $\mathrm{ch}(P)$. By Lemma \ref{lemm1}, we have $h\circ L_n(f)=\int\limits_0^1L_n(f)(x) h(x)\ dx$ converges to  $h(\lambda f)=\int\limits_{0}^1\lambda f(x)h(x)\ dx$ as $n\to\infty$ for all $h\in \mathrm{ch}(P)$ and for all $f\in L^1[0,1]$. Since every member of $K$ is an average of two members of $\mathrm{ch}(P)$, the same convergence holds for all $h\in K$. Therefore, we have
	\begin{equation}\label{lim}
		\lim\limits_{n\to\infty}\int\limits_a^bL_n(f)(x)\ dx=\int \limits_a^b \lambda f(x) \ dx,
	\end{equation}
	for all $f\in L^1[0,1]$ and $[a,b]\subset [0,1]$.\\
	\textbf{Step 2}: Let 
	\[
	G=\{g\in L^1[0,1]:g=\chi_{[a,b]}, \ [a,b]\subset [0,1]\  \textrm{or } g= \chi_{[0,1]-[a,b]}\}.
	\]
	It suffices to show that $L_n(g)$ converges to $\lambda g$ for each $g\in G$.
	If this is true, then suppose $f\in L^1[0,1]$. Since the span of the set $\{\chi_{[a,b]}:[a,b]\subset [0,1]\}$ is dense in $L^1[0,1]$, there exist a sequence of step functions $\{s_m\}$ such that 
	\[
	\lim\limits_{m\to\infty}\|s_m-f\|_1= 0.
	\]
	Thus,
	\begin{equation*}
		\begin{split}
			\|L_n(f)-\lambda f\|_1&\leq \|L_n(f)-L_n(s_m)\|_1+\|L_n(s_m)-\lambda s_m\|_1+\|\lambda s_m-\lambda f\|_1\\
			&\leq (1+\lambda)\|s_m-f\|_1+\|L_n(s_m)-\lambda s_m\|_1.
		\end{split}
	\end{equation*}
	Hence, we have $\lim\limits_{n\to\infty}\|L_n(f)-\lambda f\|_1=0$.\\
	\textbf{Step 3}: Define
	\begin{equation*}
		L_n'(g)(x)=
		\begin{cases}
			L_n(g)(x)	&\text{if } g(x)\neq 0\\
			0 &\text{if } g(x)=0.
		\end{cases}
	\end{equation*}
	Let, 
	\[
	Z=\{x\in [0,1]:g(x)=0\}.
	\]
	Then $L_n'(g)=L_n(g)\chi_{[0,1]-Z}$. 
	By Lemma \ref{lemm1}, we have, $h(L_n'(g))$ converges to $\lambda h(g)$ for each $h\in \mathrm{ch}(P)$. As observed in Step 1, every member of $K$ is an average of two members of $\mathrm{ch}(P)$, so the convergence holds for all $h\in K$ as well. We claim that
	\begin{equation}\label{equa1}
		\liminf\limits_{n\to \infty}\|L_n'(g)\|_1\geq \lambda\|g\|_1.
	\end{equation}
	Suppose this is not true, then
	\begin{equation*}
		\begin{split}
			\lim\limits_{n\to\infty}\int\limits_{0}^1L_n'(g)(x).1\ dx&=\liminf\limits_{n\to\infty}\int\limits_{0}^1L_n'(g)(x).1\ dx\\
			&\leq \liminf\limits_{n\to\infty} \|L_n'(g)\|_1<\lambda\|g\|_1=\lambda\int\limits_{0}^1g(x).1\ dx.
		\end{split}
	\end{equation*}
	Since the constant function $1\in K$, this inequality cannot arise.\\
	\textbf{Step 4:} Now, 
	\begin{equation*}
		\begin{split}
			\int\limits_0^1|L_n(g)(x)-\lambda g(x)| \ dx&=\int\limits_Z|L_n(g)(x)|\ dx+\int\limits_{[0,1]-Z}|L_n(g)(x)-\lambda|\ dx\\
			&\leq \int\limits_Z|L_n(g)(x)|\ dx+\int\limits_{[0,1]-Z}|L_n(1-g)(x)|\ dx\\
			&\quad\quad +\int\limits_{[0,1]-Z}|L_n(1)(x)-\lambda|\ dx.
		\end{split}
	\end{equation*}
	Since $\|L_n(1)-\lambda\|_1\to 0 $ as $n\to\infty$, we need to show that 
	\begin{equation}\label{eqeq}
		\lim\limits_{n\to\infty}\int\limits_Z|L_n(g)(x)|\ dx=\lim\limits_{n\to\infty}\int\limits_{[0,1]-Z}|L_n(1-g)(x)|\ dx=0.
	\end{equation}
	Suppose that $\limsup\limits_{n\to\infty}\int\limits_Z|L_n(g)(x)|\ dx>0$. Then by \ref{equa1}, $\limsup\limits_{n\to\infty}\|L_n(g)\|_1>\lambda\|g\|_1$, which is a contradiction to our assumption $1$. Therefore,  $\lim\limits_{n\to\infty}\int\limits_Z|L_n(g)(x)|\ dx=0$. Similarly, $\lim\limits_{n\to\infty}\int\limits_{[0,1]-Z}|L_n(1-g)(x)|\ dx=0$. Hence the result.
\end{proof}
\begin{theorem}\label{thm2}
	Let $\{L_n\}_{n\in\mathbb{N}}$ be a sequence of operators defined on $L^1[0,1]$ such that $\sup\limits_{n} \|L_n\|=\alpha<\infty$ and $L$ be a bounded linear operator on $L^1[0,1]$. Suppose
	\begin{enumerate}
		\item For each characteristic function $g=\chi_{[a,b]}$, where $[a,b]\subset [0,1]$, 
		\[
		\limsup\limits_{n\to\infty} \|(L_n-L+I)(g)\|_1\leq \|g\|_1,
		\]
		\item {$L_n(1)$ converges to $L(1)$,}
		\item {$L_n(p)$ converges weakly to $L(p)$ for the two functions $p=x$ and $p=x^2$}.
	\end{enumerate}
	then, $L_n(f)$ converges to $L(f)$ for all $f\in L^1[0,1]$.	
\end{theorem}
\begin{proof}
	Without loss of generality, suppose $L=I$. Otherwise, consider $\mathcal{L}_n=L_n+I-L$ and proceed.\par
	We have $\sup\limits_{n}\|L_n\|=\alpha$. We see that, whenever $\lim\limits_{n\to\infty}L_n(1)=1$, $\lim\limits_{n\to\infty}\|L_n(1)\|_1=1$, so $\sup\limits_{n}\|L_n\|\geq \sup\limits_{n}\|L_n(1)\|_1\geq 1$. Thus, we must have $\sup\limits_{n}\|L_n\|=\alpha\geq 1$. We need to prove that $\lim\limits_{n\to\infty}\|L_n(f)-f\|_1=0$ for all $f\in L^1[0, 1]$.
	
	Suppose there exists $f\in L^1[0, 1]$ such that this is not true. Then there exist a subsequence $\{L_{n_l}\}_{l\in\mathbb{N}}$ of $\{L_n\}_{n\in\mathbb{N}}$ and $\epsilon>0$ such that
	\begin{equation}\label{con1}
		\|L_{n_l}(f)-f\|_1>\epsilon
	\end{equation}
	for all $l\in\mathbb{N}$.
	Since $\sup\limits_n\|L_n\|=\alpha$, we have $\sup\limits_l\|L_{n_l}\|=\alpha_1$ for some $\alpha_1\leq\alpha$, and $\alpha_1\geq 1$, since $\sup\limits_{l}\|L_{n_l}\|\geq \sup\limits_{l}\|L_{n_l}(1)\|_1\geq 1$. Therefore, we can further find a subsequence of $\{L_{n_l}\}_{l\in\mathbb{N}}$, say $\{L_{n_k}\}_{k\in\mathbb{N}}$, such that $\lim\limits_{k\to\infty}\|L_{n_k}\|=\alpha_1$. Let $\tilde{L}_{n_k}=\dfrac{L_{n_k}}{\|L_{n_k}\|}$ and $\lambda=\dfrac{1}{\alpha_1}\leq 1$.
	
	Then, we have:
	\begin{enumerate}
		\item $\limsup\limits_{k\to\infty}\|\tilde L_{n_k}(g)\|_1\leq \lambda\|g\|_1$ for every characteristic function $g=\chi_{[a,b]}$. Indeed, by assumption (1) (with $L=I$),
		\[
		\limsup\limits_{k\to\infty}\|\tilde L_{n_k}(g)\|_1=\limsup\limits_{k\to\infty}\frac{\|L_{n_k}(g)\|_1}{\|L_{n_k}\|}\leq \frac{1}{\alpha_1}\|g\|_1=\lambda\|g\|_1,
		\]
		since $\lim\limits_{k\to\infty}\|L_{n_k}\|=\alpha_1$ and $\limsup\limits_{k\to\infty}\|L_{n_k}(g)\|_1\leq\|g\|_1$.
		
		\item $\lim\limits_{k\to\infty}\|\tilde L_{n_k}(1)-\lambda\|_1=0$. Indeed, since $\lim\limits_{k\to\infty}\|L_{n_k}(1)-1\|_1=0$ (as $L=I$, using assumption (2)) and $\lim\limits_{k\to\infty}\|L_{n_k}\|=\alpha_1$,
		\[
		\tilde L_{n_k}(1)=\frac{L_{n_k}(1)}{\|L_{n_k}\|}\longrightarrow \frac{1}{\alpha_1}=\lambda, \ as\  k\to \infty
		\]
		in $L^1$-norm.
		
		\item $\tilde L_{n_k}(p)$ converges weakly to $\lambda p$ as $k\to\infty$, for $p=x,x^2$. This follows from assumption (3) (with $L=I$), since $L_{n_k}(p)$ converges to  $p$ weakly and $\lim\limits_{k\to\infty}\|L_{n_k}\|=\alpha_1$, so for every $\psi\in (L^1[0,1])^*=L^\infty[0,1]$,
		\[
		\psi(\tilde L_{n_k}(p))=\frac{\psi(L_{n_k}(p))}{\|L_{n_k}\|}\longrightarrow \frac{\psi(p)}{\alpha_1}=\psi(\lambda p)\ as\  k\to \infty.
		\]
	\end{enumerate}
	Since $\{\tilde L_{n_k}\}_{k\in\mathbb{N}}$ is a sequence of norm-one operators satisfying the hypotheses of Lemma \ref{lemm2} with this $\lambda\in(0,1]$, we conclude
	\[
	\lim\limits_{k\to\infty}\|\tilde L_{n_k}(f)-\lambda f\|_1=0
	\]
	for all $f\in L^1[0,1]$. Multiplying by $\lim\limits_{k\to\infty}\|L_{n_k}\|=\alpha_1=1/\lambda$ and using $\tilde L_{n_k}=L_{n_k}/\|L_{n_k}\|$, this gives
	\[
	\lim\limits_{k\to\infty}\|L_{n_k}(f)-f\|_1=0
	\]
	for all $f\in L^1[0,1]$. But this contradicts \eqref{con1}. Therefore, we must have $\lim\limits_{n\to\infty}\|L_n(f)-f\|_1=0$ for all $f\in L^1[0,1]$.
\end{proof}
\begin{remark}
	Theorem~\ref{wul2} follows as a consequence of Theorem~\ref{thm2}. 
	Let $L=I$. Since $\lim\limits_{n\to\infty}\|L_n\|=1$, it follows that 
	$\sup\limits_{n}\|L_n\|<\infty$ and $\limsup\limits_{n\to\infty}\|L_n(f)\|_1\leq \|f\|_1$ 
	for all $f\in L^1[0,1]$. Thus, condition (1) is satisfied. Moreover, 
	conditions 2 and 3 of Theorem~\ref{thm2} are clearly satisfied by 
	$\{L_n\}_{n\in\mathbb{N}}$ under the hypotheses of Theorem~\ref{wul2}. 
	Hence, the result follows.
	
\end{remark}

\subsection{Grünwald's Operators and its Extension}
We consider the following sequence of operators derived from the classical literature.
\begin{definition}\label{Gn}
	Define $G_n:C[0,\pi]\rightarrow C[0,\pi]$ by 
	\[
	G_n(f)(\theta)=\frac{1}{2}\sum\limits_{k=1}^nf(\theta_k^{(n)})\{P_k(\theta-\frac{\pi}{2n})+P_k(\theta+\frac{\pi}{2n})\},
	\]
	where $P_k$ denotes the fundamental polynomials of the Lagrange interpolation operators on the Chebychev nodes $\{\cos \theta_k^{(n)}\}_{k=1}^n$, $\theta_k^{(n)}=\frac{(2k-1)\pi}{2n}$ given by,
	\[
	P_k(\theta)=(-1)^{k+1}\frac{\cos n\theta\sin \theta_k^{(n)}}{n(\cos\theta-\cos \theta_k^{(n)})}.
	\]
	for $\theta\in [0,\pi]$, $k=1,2,\dots$.
\end{definition}
Note that the function $P_k(.-\frac{\pi}{2n})+P_k(.+\frac{\pi}{2n})$ is a $2\pi$-periodic even function on $\mathbb{R}$. In $1941$, G. Grünwald proved the following result (see \cite{grunwald}).
\begin{lemma}[\textbf{G. Grünwald}, \cite{grunwald}]\label{lem}
	\[
	\frac{1}{2}\sum\limits_{k=1}^n|P_k(\theta-\frac{\pi}{2n})+P_k(\theta+\frac{\pi}{2n})|< c_1
	\]
	where $c_1>0$ is an absolute constant.
\end{lemma} 
Lemma \ref{lem} ensures the uniform boundedness of $\{G_n\}_{n\in\mathbb{N}}$. That is we have, $\|G_n(f)\|_\infty\leq c_1\|f\|_\infty$, since $c_1>0$ is independent of $n$ and $f$.
Now we prove the non-positivity of the operators $G_n$, $n=1,2,\ldots$.
\begin{theorem}
	The operators $G_n:C[0,\pi]\to C[0,\pi]$, $n=1,2\ldots$ are not positive.
\end{theorem}
\begin{proof}
	Let $\theta_k^{(n)}=\frac{(2k-1)\pi}{2n}$ for $k=1,2,\ldots n$. Then $\theta_n^{(n)}>\theta^{(n)}_{n-1}>...>\theta^{(n)}_1$. Hence $\cos\theta_n^{(n)}<\cos\theta^{(n)}_{n-1}<\ldots<\cos\theta^{(n)}_1$.\par 
	We can find a non-negative function $f\in C[0,\pi]$ such that $f(\theta^{(n)}_k)=0$ for $k=1,2,...,n-1$ and $f(\theta^{(n)}_n)>0$.
	For $\theta\in (\theta^{(n)}_{n-1}+\frac{\pi}{2n},\theta^{(n)}_n+\frac{\pi}{2n})$, we have $\theta-\frac{\pi}{2n}\in (\theta^{(n)}_{n-1},\theta^{(n)}_n)$. On this interval, the fundamental polynomial $P_n$ is negative, while $P_k(\theta-\frac{\pi}{2n})=0$ for $k=1,\ldots,n-1$. Therefore, for such $\theta$, we have 
	\[
	G_n(f)(\theta)=\frac{1}{2}f(\theta^{(n)}_n)P_n(\theta-\frac{\pi}{2n})<0.
	\]
	Thus, $G_n$ is not a positivity preserving map for $n=1,2,\ldots$ .
\end{proof}

\begin{remark}
	We note the following:
	\begin{enumerate}
		\item $\|G_n\|=\frac{1}{2}\|\sum\limits_{k=1}^n|P_k(\theta-\frac{\pi}{2n})+P_k(\theta+\frac{\pi}{2n})|\|_\infty$. To see this, we have  $\|G_n(f)\|_\infty\leq \frac{1}{2}\|\sum\limits_{k=1}^n|P_k(\theta-\frac{\pi}{2n})+P_k(\theta+\frac{\pi}{2n})\|_\infty\|f\|_\infty$. Since $P_k$'s, for $k=1,2,\ldots$ are continuous on $[0,\pi]$, there exists a $\theta_0\in[0,\pi]$ such that, $\frac{1}{2}\|\sum\limits_{k=1}^n|P_k(\theta-\frac{\pi}{2n})+P_k(\theta+\frac{\pi}{2n})|\|_\infty=\frac{1}{2}\sum\limits_{k=1}^n|P_k(\theta_0-\frac{\pi}{2n})+P_k(\theta_0+\frac{\pi}{2n})|$. Let $f\in C[0,\pi]$ be such that $\|f\|_\infty=1$ and $f(\theta_k^{(n)})=\frac{|P_k(\theta_0-\frac{\pi}{2n})+P_k(\theta_0+\frac{\pi}{2n})|}{P_k(\theta_0-\frac{\pi}{2n})+P_k(\theta_0+\frac{\pi}{2n})}$. Then, $\|G_n(f)\|_\infty=\frac{1}{2}\|\sum\limits_{k=1}^n|P_k(\theta-\frac{\pi}{2n})+P_k(\theta+\frac{\pi}{2n})|\|_\infty$.
		\item $\|G_n\|>1$. Since $P_k$'s are fundamental polynomials, we have $\frac{1}{2}\sum\limits_{k=1}^nP_k(\theta-\frac{\pi}{2n})+P_k(\theta+\frac{\pi}{2n})=1$ for all $\theta\in [0,\pi]$. Suppose $\|G_n\|=1$, then we must have $P_k(\theta+\frac{\pi}{2n})+P_k(\theta-\frac{\pi}{2n})>0$ for all $\theta\in [0,\pi]$, which is not true.
	\end{enumerate}
\end{remark} 

Now, we obtain an extension for the sequence of operators $\{G_n\}_{n\in\mathbb{N}}$ to the space $L^1(\mathbb{R})$.
Firstly, we extend the definition of $G_n$ for $n=1,2,\ldots$ to $\mathbb{R}$.
\begin{definition}\label{def1}
	Let $C(\mathbb{R})$ be the space of all continuous functions on $\mathbb{R}$ and let $f\in C(\mathbb{R})$. Suppose $j\in \mathbb{Z}$. On $[j\pi,(j+1)\pi]$, $\{\theta_k^{(n)}+j\pi\}_{k=1}^n$ is a finite sequence of nodal points. Define
	\begin{align*}
		P^j_k(\theta)&:=P_k(\theta-j\pi)=\frac{(-1)^{k+1}\cos{n(\theta-j\pi)}\sin\theta^{(n)}_k}{n(\cos(\theta-j\pi)-\cos\theta^{(n)}_k)},\\
		&\qquad\qquad\qquad\theta\in[j\pi,(j+1)\pi],\quad j\in\mathbb Z,
	\end{align*}
	and
	\begin{align*}
		G^j_n(f)(\theta)&:=\frac12\sum_{k=1}^nf(\theta^{(n)}_k+j\pi)\left(P^j_k(\theta+\frac\pi{2n})+P^j_k(\theta-\frac\pi{2n})\right),\\
		&\qquad\qquad\qquad\theta\in[j\pi,(j+1)\pi),\quad j\in\mathbb Z.
	\end{align*}
	For $f\in C(\mathbb{R})$, we define
	\[
	G_n(f)(\theta)=G_n^j(f) (\theta)
	\]
	for $j\in \mathbb{N}$ and $\theta\in [j\pi,(j+1)\pi)$. 
\end{definition}
\begin{remark}
	Note that we use the same notation $G_n(f)$ as in Definition \ref{Gn}. $G_n(f)$ is defined almost everywhere on $\mathbb{R}$. 
\end{remark}
\begin{definition}\label{oper}
	For $f\in L^1(\mathbb{R})$, let $\widehat{f}$ denote the Fourier transform of $f$ and $\widecheck{f}$, the inverse Fourier transform of $f$. It is also well known that for each $f\in L^1(\mathbb{R})$, $\hat{f}$ is in $C(\mathbb{R})$. Let
	\begin{multline*}
		U:=\{f\in L^1(\mathbb{R}):\widehat{f}\in L^1(\mathbb{R})\ \textrm{and } \sum_{l=-\infty}^{\infty}|\widehat{f}(d_l)|<\infty\ \textrm{where }\\
		|\widehat{f}(d_l)|=\sup\limits_{\theta\in [l\pi,(l+1)\pi]}|\widehat{f}(\theta)|\},
	\end{multline*}
	where $d_l$ denotes the point where $|\hat{f}|$ attains its maximum in the interval $[l\pi,(l+1)\pi]$.
	For $f\in U$, define
	\[
	K_n(f)=\widecheck{G_n\widehat{(f)}}.
	\]
\end{definition}
\begin{remark}
	The set $U$ is a technical condition ensuring that the operator $K_n$ is well-defined on a sufficiently large subset of $L^1(\mathbb{R})$ where the Fourier transform can be manipulated conveniently.
\end{remark}
We see that $K_n(f)$ is locally integrable for all $n \in \mathbb{N}$ and $f \in U$. Thus, $K_n$ maps each $f \in U$ to the space $L_{loc}^1(\mathbb{R})$.

We extend $K_n$, for $n=1,2,\ldots$ to $L^1(\mathbb{R})$ employing the standard techniques of using an approximate identity.
\begin{definition}
	An approximate identity is a family of functions $\phi_\delta \in L^1(\mathbb{R})$, $\delta>0$ such that
	\begin{enumerate}
		\item The family is bounded in the $L^1$-norm, that is, there is a $C > 0$ such that $\|\phi_\delta\|_{L^1(\mathbb{R})} \leq C$, for all $\delta > 0$.
		\item $\int_{\mathbb{R}} \phi_\delta(x) \, dx = 1$, for all $\delta > 0$.
		\item For any $\epsilon > 0$, $\lim_{\delta \to 0} \int_{|x| \geq \epsilon} |\phi_\delta(x)| \, dx = 0$.
	\end{enumerate}
\end{definition}

Suppose $\{\phi_\delta\}_{\delta>0}$ be an approximate identity such that $\hat{\phi_\delta}$ is compactly supported. Then, for all $f\in L^1(\mathbb{R})$, $\widehat{f*\phi_{\delta}}=\widehat{f}\widehat{\phi_\delta}$ is also compactly supported (here $*$ denotes the operation convolution).

\begin{definition}\label{hndelta}
	Suppose $f\in L^1(\mathbb{R})$. Let $\{\phi_\delta\}_{\delta>0}$ be an approximate identity such that $\widehat{\phi_\delta}$ is compactly supported for each $\delta>0$.
	Define 
	\[
	H_{n,\delta}(f)(x)=\frac{1}{2\pi}\sum\limits_{l=-\infty}^\infty\frac{1}{2}\sum\limits_{k=1}^n\widehat{f*\phi_\delta}(\theta_k^{(n)}+l\pi)\int\limits_{l\pi}^{(l+1)\pi}(P_k^l(\theta+\frac{\pi}{2n})+P_k^l(\theta-\frac{\pi}{2n}))e^{ix\theta}d\theta.
	\]	
\end{definition}
\begin{remark}\label{rem}
	Here $H_{n,\delta}$ is well defined. Since $\widehat{f*\phi_\delta}$ is compactly supported the series appearing in the Definition \ref{hndelta} is actually a finite sum. We also observe that 
	\[
	H_{n,\delta}(f)=K_n(f*\phi_\delta)
	\]
	for $f\in L^1(\mathbb{R})$.
\end{remark}

\subsubsection{Convergence of $\{G_n\}_{n\in\mathbb{N}}$ and $\{H_{n,\delta}\}_{n\in\mathbb{N},\delta>0}$}

In this subsection, apply Corollary \ref{cor5.10} on $\{G_n\}_{n\in\mathbb{N}}$ and as an application we prove the convergence of $\{H_{n,\delta}\}_{n\in\mathbb{N},\delta>0}$. First, we apply Corollary \ref{cor5.10} to $\{G_n\}_{n\in\mathbb{N}}$.
\begin{proposition}\label{corol1}
	The sequence of operators $\{G_n\}_{n\in\mathbb{N}}$ satisfies the following convergence:
	\[
	\lim\limits_{n\to\infty}\|G_n(f)-f\|_\infty=0
	\]
	for all $f\in C[0,\pi]$
\end{proposition}
\begin{proof}
	We already have $\sup\limits_n \|G_n\|=\alpha$ (Lemma \ref{lem}), where $\alpha>1$ (since $G_n(1)=1$ for all $n\in \mathbb{N}$, $\|G_n\|>1$). Now we apply Corollary \ref{cor5.10}. Here $E=C[0,\pi]$ and $D=\{p(\cos):p\ \text{is a polynomial on}\ [-1,1]\}$, which is dense in $C[0,\pi]$. First, we show that 
	\begin{equation}\label{limsup}
		\limsup\limits_{n\to\infty}\|G_n(f)\|_{\infty}\leq \|f\|_\infty,
	\end{equation}
	for all $f\in D$. Let $\theta\in [0,\pi]$ and $m\in \mathbb{N}$. We have 
	\[
	G_n(p(\cos))(\theta)=\frac{1}{2}\sum\limits_{k=1}^n (p(\cos\theta_k^{(n)})-p(\cos\theta))(P_k(\theta+\frac{\pi}{2n})+P_k(\theta-\frac{\pi}{2n}))+p(\cos\theta).
	\]
	We have $|p(\cos \theta)-p(\cos\eta)|\leq \sup\limits_{x\in [-1,1]}|p'(x)||\cos\theta-\cos\eta|=K|\cos\theta-\cos\eta|$, for all $\theta,\eta\in [0,\pi]$.
	Thus, we have 
	\begin{align*}
		|G_n(p(\cos&))(\theta)|\\
		&\leq\frac{1}{2}\sum\limits_{k=1}^n |p(\cos\theta_k^{(n)})-p(\cos\theta)||P_k(\theta+\frac{\pi}{2n})+P_k(\theta-\frac{\pi}{2n})|+|p(\cos\theta)|\\
		&\leq \frac{K}{2}\sum\limits_{k=1}^n |\cos\theta_k^{(n)}-\cos\theta||P_k(\theta+\frac{\pi}{2n})+P_k(\theta-\frac{\pi}{2n})|+|p(\cos\theta)|\\
		&\leq \frac{K}{2}\sup\limits_{\theta\in [0,\pi]}\sum\limits_{k=1}^n |\cos\theta_k^{(n)}-\cos\theta||P_k(\theta+\frac{\pi}{2n})+P_k(\theta-\frac{\pi}{2n})|+\|p(\cos)\|_\infty.
	\end{align*}
	Hence, the condition is verified. We need to only show that 
	\[
	\frac{1}{2}\sum\limits_{k=1}^n |\cos\theta_k^{(n)}-\cos\theta||P_k(\theta+\frac{\pi}{2n})+P_k(\theta-\frac{\pi}{2n})|
	\]
	converges to $0$ uniformly on $[0,\pi]$.
	Let $\epsilon>0$ and $\theta\in [0,\pi]$. Consider the two intervals $|\theta-\theta_k^{(n)}|>\epsilon$ and $|\theta-\theta_k^{(n)}|\leq \epsilon$. 
	For $k$ such that $|\theta-\theta_k^{(n)}|>\epsilon$, we have the following inequality from \cite{grunwald}
	\[
	\frac{1}{2}|P_k(\theta+\frac{\pi}{2n})+P_k(\theta-\frac{\pi}{2n})|\leq \frac{\pi^3}{4n^2}\frac{1}{(\theta-\theta_k^{(n)}-\frac{\pi}{2n})^2}.
	\]
	Therefore, we have  $\frac{1}{2}|P_k(\theta+\frac{\pi}{2n})+P_k(\theta-\frac{\pi}{2n})|=\mathcal{O}(\frac{1}{n^2})$. Thus, there exists an $n_0\in\mathbb{N}$ such that
	\[
	\frac{1}{2}\sum\limits_{|\theta-\theta_k^{(n)}|>\epsilon}|P_k(\theta+\frac{\pi}{2n})+P_k(\theta-\frac{\pi}{2n})|\leq \frac{\epsilon}{2},
	\]
	for all $n\geq n_0$. 
	Now
	\begin{align*}
		\frac{1}{2}\sum\limits_{k=1}^n &|\cos\theta_k^{(n)}-\cos\theta||P_k(\theta+\frac{\pi}{2n})+P_k(\theta-\frac{\pi}{2n})|\\
		&=\frac{1}{2}\left (\sum\limits_{|\theta-\theta_k^{(n)}|\leq\epsilon}+\sum\limits_{|\theta-\theta_k^{(n)}|>\epsilon}\right ) |\cos\theta_k^{(n)}-\cos\theta||P_k(\theta+\frac{\pi}{2n})+P_k(\theta-\frac{\pi}{2n})|\\
		&\leq \alpha\epsilon+\epsilon.
	\end{align*}
	for all $n\geq n_0$, since we have $|\cos\theta_k^{(n)}-\cos\theta|\leq 2$ for all $k=1,2,\ldots,n$  for $|\theta_k^{(n)}-\theta|>\epsilon$ and $|\cos\theta_k^{(n)}-\cos\theta|\leq|\theta_k^{(n)}-\theta|\leq \epsilon$, for $|\theta_k^{(n)}-\theta|\leq \epsilon$. We have 
	\[
	\|G_n(p(\cos))\|_\infty\leq K(\alpha+1)\epsilon+\|p(\cos)\|_\infty.
	\]
	for all $n\geq n_0$. Hence (\ref{limsup}) is verified.
	
	Thus, it suffices to prove the convergence only on the test set $\{1,\cos, \cos^2\}$. But this is obvious from the above observations.
	
\end{proof}
\begin{remark}
	This result can be extended to the sequence of operators $\{G_n^l\}_{n\in\mathbb{N}}$ defined on $[l\pi,(l+1)\pi]$ for each $l\in \mathbb{Z}$ in a similar way.
\end{remark}
We recall the family of operators $H_{n,\delta}:L^1(\mathbb{R})\to L^1(\mathbb{R})$. For $f\in L^1[0,\pi]$, then we have $f\chi_{[0,\pi]}\in L^1(\mathbb{R})$, so that we can restrict the definition to the space $L^1[0,\pi]$. We have the following convergence result.
\begin{proposition}
	Consider the family of operators $H_{n,\delta}:L^1[0,\pi]\to L^1[0,\pi]$. Then, we have 
	\[
	\lim\limits_{\delta\to 0}\lim\limits_{n\to\infty}\|H_{n,\delta}(f)-f\|_1=0.
	\]
\end{proposition}
\begin{proof}
	
	For $x\in [0,\pi]$, we have,
	\[
	|H_{n,\delta}(f)(x)-f*\phi_\delta(x)|\leq \frac{1}{2\pi}\sum\limits_{l=-m_\delta}^{m_\delta}\|G_n^l(\widehat{f*\phi_\delta})-\widehat{f*\phi_\delta}\|_\infty
	\]
	Hence,
	\[
	\|H_{n,\delta}(f)-f*\phi_\delta\|_1\leq \frac{1}{2}\sum\limits_{l=-m_\delta}^{m_\delta}\|G_n^l(\widehat{f*\phi_\delta})-\widehat{f*\phi_\delta}\|_\infty.
	\]
	By Proposition \ref{corol1}, The right hand side tends to $0$ as $n\to \infty$ for all $f\in L^1[0,\pi]$. Now, 
	$\|H_{n,\delta}(f)-f\|_1\leq \|H_{n,\delta}(f)-f*\phi_\delta\|_1+\|f*\phi_\delta-f\|_1$. Thus, we have 
	$\|H_{n,\delta}(f)-f\|_1$ converges to $0$ as $n\to\infty$ and $\delta\to0$.
\end{proof}

\section{A Direct Approach for the convergence of $\{G_n\}$ and $\{H_{n,\delta}\}$}\label{direct}
In $1941$, G. Grünwald proved the following convergence result.
\begin{theorem}[\textbf{G. Grünwald}, \cite{grunwald}]\label{grun}
	Let $f\in C[-1,1]$, then 
	\[
	\lim\limits_{n\rightarrow\infty}\frac{1}{2}\{L_n(f)(\theta-\frac{\pi}{2n})+L_n(f)(\theta+\frac{\pi}{2n})\}=f(\cos\theta),
	\]
	uniformly on $[0,\pi]$.
\end{theorem}
Here $L_n$ denotes the Lagrange interpolation operator on the Chebychev nodes. This can be formulated in terms of $G_n$ as follows:
For $f\in C[-1,1]$,
\[
\lim\limits_{n\to\infty}G_n(f(\cos))=f(\cos)
\]
uniformly on $[0,\pi]$. We also observe that the proof of Theorem \ref{grun} equally works for $\{G_n\}_{n\in\mathbb{N}}$, that is, we have,
$\lim\limits_{n\to\infty}G_n(f)=f$ in $C[0,\pi]$.\par

We recall the definition of modulus of continuity below. 
\begin{definition}
	Let $f:I\to \mathbb{R}$ be a bounded functions on a real interval $I$. For $\delta>0$,  the modulus of continuity $\omega(f,\delta)$ of $f$ with argument $\delta>0$ is defined as $
	\omega(f,\delta):=sup\{|f(x)-f(y)|: x,y\in I,\ |x-y|\leq\delta\}.$
\end{definition}

The following quantitative estimate was obtained for the convergence of the Gr\"unwald operators in \cite{mills1}.

\begin{theorem}\cite{mills1}\label{moc}
	Let $f\in C[-1,1]$ and $x\in [-1,1]$. Then, we have
	\[
	\left|G_n(f\circ\cos)(\theta)-f(\cos\theta)\right|
	=
	\mathcal{O}\left(
	\omega\left(f,\frac{\sqrt{1-x^2}}{n}\right)
	+
	\omega\left(f,\frac{1}{n^2}\right)
	\right),
	\]
	where $x=\cos\theta$.
\end{theorem}

In the following theorem, we obtain convergence of $K_n(f)$ to $f$ on compact intervals of $\mathbb{R}$ in the $L^1$-norm. We emphasize that the majority of information concerning a function in $L^1(\mathbb{R})$ is concentrated within a large compact set, with the norm outside this compact interval being negligible. Consequently, convergence within compact subsets of $\mathbb{R}$ holds much significance.\par 	
\begin{theorem}\label{extn}
	Let $U$ and $K_n$ be as in Definition \ref{oper} and suppose $f\in U$. Then
	\[
	\lim\limits_{n\rightarrow\infty}\|(K_n(f)-f)\chi_{[a,b]}\|_1=0
	\]
	for any compact interval $[a,b]$ of $\mathbb{R}$. 
\end{theorem}

\begin{proof}
	Let $f\in U$.	We see that
	\begin{equation*}
		\begin{split}
			\|G_n(\widehat{f})\|_1
			&=\lim\limits_{m\rightarrow \infty}\sum\limits_{l=-m}^m\int\limits_{l\pi}^{(l+1)\pi}|G_n^l(\widehat{f})(\theta)|\ d\theta\\
			&\leq\lim\limits_{m\rightarrow \infty}\sum\limits_{l=-m}^m\frac{1}{2}\sum_{k=1}^n|\widehat{f}(\theta_k^{(n)}+l\pi)|\int\limits_{l\pi}^{(l+1)\pi}|P_k^l(\theta-\frac{\pi}{2n})+P_k^l(\theta+\frac{\pi}{2n})|\ d\theta\\
			&\leq \lim\limits_{m\rightarrow \infty}\sum\limits_{l=-m}^m|\widehat{f}(d_l)|\int\limits_{l\pi}^{(l+1)\pi}\frac{1}{2}\sum_{k=1}^n|P_k^l(\theta-\frac{\pi}{2n})+P_k^l(\theta+\frac{\pi}{2n})|\ d\theta.
		\end{split}
	\end{equation*} 
	We observe that the integrand in the above term is bounded by the same absolute constant $c_1$ for each $l\in\mathbb{Z}$, by Lemma \ref{lem}. Since $f\in U$, we have that $G_n(\widehat{f})\in L^1(\mathbb{R})$. Hence the inverse Fourier transform of $G_n(\widehat{f})$ exists. Now,
	\begin{equation*}
		\begin{split}
			\widecheck{G_n\widehat{(f)}}(x)-f(x)&=\frac{1}{2\pi}\int\limits_{-\infty}^{\infty}G_n(\widehat{f})(\theta)e^{ix\theta}d\theta-f(x)\\
			&=\frac{1}{2\pi}\int\limits_{-\infty}^{\infty}(G_n(\widehat{f})(\theta)-\widehat{f}(\theta))e^{ix\theta}d\theta
		\end{split}
	\end{equation*}
	almost everywhere. We use the fact that $f(x)= \int\limits_{-\infty}^{\infty}\widehat{f}(\theta)e^{ix\theta}d\theta$ almost everywhere by the Fourier inversion theorem.
	We have,
	\begin{equation}\label{eqn}
		|\widecheck{G_n\widehat{(f)}}(x)-f(x)|\leq \frac{1}{2\pi}\int\limits_{-\infty}^{\infty}|G_n(\widehat{f})(\theta)-\widehat{f}(\theta)|\ d\theta
	\end{equation}
	almost everywhere on $\mathbb{R}$.
	
	By our assumption both $G_n(\widehat{f})$ and $\widehat{f}$ are in $L^1(\mathbb{R})$. We observe that $|G_n(\widehat{f})(\theta)|\leq c_1\sum\limits_{l=-\infty}^{\infty}|\widehat{f}(d_l)|\chi_{[l\pi,(l+1)\pi]}(\theta)$ for $\theta\in \mathbb{R}$. The function in the right hand side is an integrable function over $\mathbb{R}$ by our assumption. Moreover, by Grünwald's theorem, $G_n(\widehat{f})(\theta)\to \widehat{f}(\theta)$ for all $\theta\in\mathbb{R}$. Since $\widehat{f}\in L^1(\mathbb{R})$, the integrand in \eqref{eqn} is dominated by an $L^1$ function independent of $n$. Therefore, by the Lebesgue Dominated Convergence Theorem, the right-hand side of \eqref{eqn} converges to $0$ as $n\to\infty$. This implies $K_n(f)(x)\to f(x)$ almost everywhere on $\mathbb{R}$. Hence the result follows.
\end{proof}
\begin{corollary}
	For $f\in L^1(\mathbb{R})$, if $\widehat{f}$ is a compactly supported function, then we have
	\[
	\lim\limits_{n\rightarrow\infty}\|(K_n(f)-f)\chi_{[a,b]}\|_1=0
	\]
	for any compact interval $[a,b]$ of $\mathbb{R}$.  
\end{corollary}	

\begin{proof}
	Since $\widehat{f}$ is compactly supported, the series in the definition of $U$ reduces to a finite sum and hence $f\in U$. The result follows from Theorem \ref{extn}.
\end{proof}
\begin{corollary}\label{cor}
	Let $f\in C^2(\mathbb{R})\cap L^1(\mathbb{R})$ and $f',f''\in L^1(\mathbb{R})$. Then ${f}\in U$ and thus 
	\[
	\lim\limits_{n\rightarrow\infty}\|(K_n(f)-f)\chi_{[a,b]}\|_1=0
	\]
	for any compact interval $[a,b]$ of $\mathbb{R}$.  
\end{corollary}
\begin{proof}
	Note that whenever $f\in C^2(\mathbb{R})\cap L^1(\mathbb{R})$ and $f',f''\in L^1(\mathbb{R})$, we have Fourier transform of $f$, $\widehat{f}$ exists. We know that for $x \in \mathbb{R}$
	\[
	\widehat{f''}(x)=-x^2.\widehat{f}(x).
	\]
	Let $g(x)=f(x)-f''(x)$. Then $\widehat{g}(x)=\widehat{f}(x)-\widehat{f''}(x)$ and 
	$\widehat{g}(x)=(1+x^2)\widehat{f}(x)$, so that $|\widehat{f}(x)|\leq \frac{1}{1+x^2}|\widehat{g}(x)|$. In this inequality, the right hand side is the product of two functions, 
	$\frac{1}{1+x^2}$, which is an integrable function and $\widehat{g}\in C_0(\mathbb{R})$ by the Riemann Lebesgue Lemma(page $249$, \cite{folland}). Here, the right hand side is in $L^1(\mathbb{R})$ hence the left hand side. Therefore, $\widehat{f}\in L^1(\mathbb{R})$. Again, 
	\[
	\sum\limits_{l=-\infty}^{\infty}|\widehat{f}(d_l)|\leq C \sum\limits_{l=-\infty}^{\infty} \frac{1}{1+d_l^2},
	\]
	where $C>0$ is a constant. Note that each positive integer $l$, $d_l\in [l\pi,(l+1)\pi]$, so that $d_l\geq l\pi$ and $\frac{1}{1+d_l^2}\leq \frac{1}{1+(l\pi)^2} \leq \frac{1}{1+l^2}$. For negative integers, we have 
	$d_{-l}\leq (-l+1)\pi$ and $-d_{-l}\geq (l-1)\pi$ where $d_{-l}\in [-l\pi,(-l+1)\pi]$. So that $\frac{1}{1+d_{-l}^2}\leq \frac{1}{1+(l-1)^2}$. Thus the series in the right hand side converges.
\end{proof}	
\begin{remark}
	If $f$ is compactly supported then all its derivatives are compactly supported and hence they belong to $L^1(\mathbb{R})$. Let $f$ be compactly supported, twice differentiable such that $f''\in C(\mathbb{R})$. From Corollary \ref{cor}, we see that $f\in U$ and hence collection of all such functions are in $U$. Moreover, this collection is dense in $L^1(\mathbb{R})$. Thus Grünwald-type theorem holds in a dense class of $L^1(\mathbb{R})$.\par
	We also observe that the set $\tilde{U}=\{f\in L^1(\mathbb{R}):\widehat{f}\ \text{is compactly supported}\}$ is dense in $L^1(\mathbb{R})$. \par
	Note that when $f\in C(\mathbb{R})\cap U$, then we have $G_n(f)(\theta)\rightarrow f(\theta)$ almost every $\theta\in \mathbb{R}$. Also $K_n(f)(\theta)=\widecheck{G_n\widehat{(f)}} 
	(\theta)\rightarrow f(\theta)$ almost everywhere on $\mathbb{R}$. Therefore, $G_n(f)$ and $K_n(f)$ are asymptotically equivalent in $L^1$ sense on compact intervals of $\mathbb{R}$.
\end{remark}
\begin{remark}
	The operators $K_n$, $n=1,2,\ldots$ are defined in terms of the Fourier transform of $f$ on some grid points in $\mathbb{R}$. Hence if for $f\in U$, $\hat{f}$ is known at some grid points, then $f$ can be approximated on all compact subsets of $\mathbb{R}$ in the $L^1$ sense.
\end{remark}
\subsection{Computations}
In this subsection, we compute the integrals occurring in the extended operator. Thus we get an explicit formula for our new operator.\par
Let $f\in U$ and $p\in \mathbb{R}$,
then,
\[ K_n(f)(p)=\frac{1}{2\pi}\sum\limits_{l=-\infty}^{\infty}\frac{1}{2}\sum_{k=1}^{n}\widehat{f}(\theta_k^{(n)}+l\pi)\int\limits_{l\pi}^{(l+1)\pi}(P_k^l(\theta+\frac{\pi}{2n})+P_k^l(\theta-\frac{\pi}{2n}))e^{ip\theta}\ d\theta.
\]
We compute the integral in this definition for each $l\in \mathbb{Z}$. Let $\theta\in [l\pi, (l+1)\pi]$. We have
\begin{equation*}
	P_k^l(\theta+\frac{\pi}{2n})+P_k^l(\theta-\frac{\pi}{2n})
	=P_k(\theta+\frac{\pi}{2n}-l\pi)+P_k(\theta-\frac{\pi}{2n}-l\pi)
\end{equation*}
Note that for $\theta\in [l\pi,(l+1)\pi]$, $\eta=\theta-l\pi\in [0,\pi]$.
We evaluate the integral
\[
\int\limits_{l\pi}^{(l+1)\pi}(P_k^l(\theta-\frac{\pi}{2n})+P_k^l(\theta+\frac{\pi}{2n}))e^{ip\theta}  \ d\theta=(-1)^l\int\limits_{0}^{\pi}(P_k(\eta-\frac{\pi}{2n})+P_k(\eta+\frac{\pi}{2n}))e^{ip\eta}\  d\eta,
\]
on applying a change of variable, $\eta=\theta-l\pi$.

We evaluate the integrals\\ $\displaystyle{\int\limits_{0}^{\pi}(P_k(\eta-\frac{\pi}{2n})+P_k(\eta+\frac{\pi}{2n}))\cos (p\eta) \ d\eta}$ and $\displaystyle{\int\limits_{0}^{\pi}(P_k(\eta-\frac{\pi}{2n})+ P_k(\eta+\frac{\pi}{2n})) \sin (p\eta) \ d\eta}$. We note that 
We have
\begin{multline*}
	P_k(\theta) = \frac{(-1)^{k+1} \sin \theta_k^{(n)}}{n} \cdot (\cos\theta - \cos \theta_1^{(n)}) \cdots (\cos\theta - \cos \theta_{k-1}^{(n)}) \\
	\times (\cos\theta - \cos \theta_{k+1}^{(n)}) \cdots (\cos\theta - \cos \theta_n^{(n)}).
\end{multline*}

The product simplifies to a polynomial in $\cos\theta$:
\begin{multline*}
	(\cos\theta - \cos \theta_1^{(n)}) \cdots (\cos\theta - \cos \theta_{k-1}^{(n)})(\cos\theta - \cos \theta_{k+1}^{(n)}) \cdots (\cos\theta - \cos \theta_n^{(n)}) = \\
	\cos^{n-1}\theta - \left(\sum_{\substack{j = 1\\ j \ne k}}^n \cos \theta_j^{(n)}\right) \cos^{n-2}\theta + \cdots + (-1)^{n-1} \prod_{\substack{j = 1\\ j \ne k}}^n \cos \theta_j^{(n)}.
\end{multline*}

Let \( \tilde{\theta} = \theta - \frac{\pi}{2n} \). To evaluate
\[
\int_0^{\pi} P_k(\tilde{\theta}) \cos(p\theta) \, d\theta,
\]
it suffices to compute the integrals
\[
\int_0^{\pi} \cos^r \tilde{\theta} \cos(p\theta) \, d\theta \quad \text{for } r = 0,1,\dots,n-1.
\]

We proceed by repeated use of the identity: $2\cos A \cos B = \cos(A + B) + \cos(A - B).$ Thus, for \( r \ge 1 \),
\begin{align*}
	\int_0^{\pi} \cos^r \tilde{\theta} \cos(p\theta) \, d\theta 
	&= \frac{1}{2} \int_0^{\pi} \cos^{r-1} \tilde{\theta} \cdot 2 \cos \tilde{\theta} \cos(p\theta) \, d\theta \\
	&= \frac{1}{2} \int_0^{\pi} \cos^{r-1} \tilde{\theta} \left[ \cos(\tilde{\theta} + p\theta) + \cos(\tilde{\theta} - p\theta) \right] d\theta \\
	&= \frac{1}{2^2} \int_0^{\pi} \cos^{r-2} \tilde{\theta} \cos\left( \frac{3\tilde{\theta} + p\theta}{2} \right) d\theta \\
	&\quad + \frac{1}{2^2} \int_0^{\pi} \cos^{r-2} \tilde{\theta} \cos\left( \frac{\tilde{\theta} + p\theta}{2} \right) d\theta \\
	&\quad + \frac{1}{2^2} \int_0^{\pi} \cos^{r-2} \tilde{\theta} \cos\left( \frac{3\tilde{\theta} - p\theta}{2} \right) d\theta \\
	&\quad + \frac{1}{2^2} \int_0^{\pi} \cos^{r-2} \tilde{\theta} \cos\left( \frac{\tilde{\theta} - p\theta}{2} \right) d\theta \\
	&\quad \vdots
\end{align*}

Each iteration reduces the power of \( \cos \tilde{\theta} \) and breaks the integrand into a sum of cosines with linear combinations of \( \tilde{\theta} \) and \( p\theta \). Ultimately, the integral reduces to sums of integrals of cosine functions.

Continuing this recursive process, we have \( 2^r \) such integrals at the \( r \)th step. Let
\[
S_r = \left\{ \phi \ \middle| \ \phi : \{1, 2, \dots, r\} \to \{-1,1\} \right\},
\]
so that \( |S_r| = 2^r \) for each \( r = 1, \dots, n-1 \).

Define
\[
Q(\phi, r) = 2^{r-1} + \phi(1)\cdot 2^{r-2} + \phi(2)\cdot 2^{r-3} + \cdots + \phi(r-1)\cdot 1.
\]

Then,
\begin{align*}
	\int_0^\pi \cos^r \tilde{\theta} \cos(p\theta)\, d\theta
	&= \frac{1}{2^r} \sum_{\phi \in S_r} \int_0^\pi \cos\left( \frac{(Q(\phi,r) + p\phi(r))\theta}{2^r} - \frac{\pi}{2n} \cdot \frac{Q(\phi,r)}{2^r} \right) d\theta \\
	&= \sum_{\phi \in S_r} \frac{1}{Q(\phi,r) + p\phi(r)} \left\{
	\sin\left( \frac{(Q(\phi,r) + p\phi(r))\pi}{2^r} - \frac{\pi}{2n} \cdot \frac{Q(\phi,r)}{2^r} \right) \right. \\
	&\quad\quad\quad\quad + \left. \sin\left( \frac{\pi}{2n} \cdot \frac{Q(\phi,r)}{2^r} \right) \right\}.
\end{align*}

Now, let \( \theta' = \theta + \frac{\pi}{2n} \). Then, for \( r = 1,2,\ldots,n-1 \), we get
\begin{align*}
	\int_0^\pi \cos^r \theta' \cos(p\theta)\, d\theta
	&= \frac{1}{2^r} \sum_{\phi \in S_r} \int_0^\pi \cos\left( \frac{(Q(\phi,r) + p\phi(r))\theta}{2^r} + \frac{\pi}{2n} \cdot \frac{Q(\phi,r)}{2^r} \right) d\theta \\
	&= \sum_{\phi \in S_r} \frac{1}{Q(\phi,r) + p\phi(r)} \left\{
	\sin\left( \frac{(Q(\phi,r) + p\phi(r))\pi}{2^r} + \frac{\pi}{2n} \cdot \frac{Q(\phi,r)}{2^r} \right) \right. \\
	&\quad\quad\quad\quad - \left. \sin\left( \frac{\pi}{2n} \cdot \frac{Q(\phi,r)}{2^r} \right) \right\}.
\end{align*}

Thus, the sum of these integrals is
\begin{multline*}
	\int_0^\pi \left( \cos^r\left(\theta - \frac{\pi}{2n}\right) + \cos^r\left(\theta + \frac{\pi}{2n}\right) \right) \cos(p\theta)\, d\theta \\
	= 2 \sum_{\phi \in S_r} \frac{1}{Q(\phi,r) + p\phi(r)} \sin\left( \frac{(Q(\phi,r) + p\phi(r))\pi}{2^r} \right) \cos\left( \frac{\pi}{2n} \cdot \frac{Q(\phi,r)}{2^r} \right).
\end{multline*}

For \( r = 0 \),
\[
\int_0^\pi \left( \cos^0\left(\theta - \frac{\pi}{2n}\right) + \cos^0\left(\theta + \frac{\pi}{2n}\right) \right) \cos(p\theta)\, d\theta =
\begin{cases}
	2\pi & \text{if } p = 0, \\
	\displaystyle \frac{2 \sin(p\pi)}{p} & \text{if } p \neq 0.
\end{cases}
\]

Now we evaluate \( \displaystyle \int_0^\pi P_k(\tilde{\theta}) \sin(p\theta)\, d\theta \).

For \( r = 1,2,\ldots,n-1 \),
\begin{multline*}
	\int_0^\pi \left( \cos^r\left(\theta - \frac{\pi}{2n}\right) + \cos^r\left(\theta + \frac{\pi}{2n}\right) \right) \sin(p\theta)\, d\theta \\
	= 2 \sum_{\phi \in S_r} \frac{1}{Q(\phi,r) + p\phi(r)} \cos\left( \frac{\pi}{2n} \cdot \frac{Q(\phi,r)}{2^r} \right)
	\left[ 1 - \cos\left( \frac{(Q(\phi,r) + p\phi(r))\pi}{2^r} \right) \right].
\end{multline*}

For \( r = 0 \),
\[
\int_0^\pi \left( \cos^0\left(\theta - \frac{\pi}{2n}\right) + \cos^0\left(\theta + \frac{\pi}{2n}\right) \right) \sin(p\theta)\, d\theta =
\begin{cases}
	0 & \text{if } p = 0, \\
	\displaystyle \frac{2(1 - \cos p\pi)}{p} & \text{if } p \neq 0.
\end{cases}
\]
This subsection can be summarized as follows:\\

Let $f \in U$ and $p \in \mathbb{R}$. Then,
\[
K_n(f)(p) = \frac{1}{2\pi} \sum_{l=-\infty}^{\infty} \frac{1}{2}(-1)^l \sum_{k=1}^n \widehat{f}(\theta_k^{(n)} + l\pi) V_k(p),
\]
where
\[
V_k(p) = \int_0^\pi \left(P_k\left(\eta - \frac{\pi}{2n}\right) + P_k\left(\eta + \frac{\pi}{2n}\right)\right) e^{ip\eta} \, d\eta = \sum_{r=0}^{n-1} s_r^k R_r(p),
\]
with the coefficients $s_r^k$ defined by:
\begin{align*}
	s_0^k &= 1, \\
	s_1^k &= -\sum_{\substack{1 \leq j \leq n \\ j \ne k}} \cos \theta_j^{(n)}, \\
	s_2^k &= \sum_{\substack{1 \leq j_1, j_2 \leq n \\ j_1 \ne j_2,\, j_1,j_2 \ne k}} \cos \theta_{j_1}^{(n)} \cos \theta_{j_2}^{(n)}, \\
	&\vdots \\
	s_{n-1}^k &= (-1)^{n-1} \prod_{\substack{1 \leq j \leq n \\ j \ne k}} \cos \theta_j^{(n)}.
\end{align*}

Each $R_r(p)$ is given by
\[
R_r(p) = R_r^1(p) + i R_r^2(p),
\]
where
\[
R_0(p) = 
\begin{cases}
	2 & \text{if } p = 0, \\
	\frac{2 \sin p\pi}{p} + i \frac{2(1 - \cos p\pi)}{p} & \text{otherwise},
\end{cases}
\]
and for $r = 1,2,\dots, n-1$,
\begin{align*}
	R_r^1(p) &= 2 \sum_{\phi \in S_r} \frac{1}{Q(\phi, r) + p \phi(r)} 
	\sin\left( \frac{(Q(\phi, r) + p \phi(r)) \pi}{2^r} \right)
	\cos\left( \frac{\pi}{2n} \cdot \frac{Q(\phi, r)}{2^r} \right), \\
	R_r^2(p) &= 2 \sum_{\phi \in S_r} \frac{1}{Q(\phi, r) + p \phi(r)} 
	\cos\left( \frac{\pi}{2n} \cdot \frac{Q(\phi, r)}{2^r} \right)
	\left(1 - \cos\left( \frac{(Q(\phi, r) + p \phi(r)) \pi}{2^r} \right)\right).
\end{align*}

Here,
\[
S_r = \left\{ \phi \mid \phi: \{1,2,\dots,r\} \to \{-1,1\} \right\}, \quad \text{and} \quad 
Q(\phi, r) = 2^{r-1} + \sum_{j=1}^{r-1} \phi(j) \cdot 2^{r-1-j}.
\]
\par
In the following part, we give the numerically computed values of $K_{n,m}(f)(p)$, the partial sums of the series in the expression of $K_n(f)(p)$. We compute the numerical values for given $n,m$ and $p$ values and $f(x)=e^{-x^2}$.
\begin{center}
	\begin{tabular}{ |p{2cm}|p{2cm}|p{1.5cm}|p{6 cm}| }
		\hline
		&&&\\
		$n$& $m$& p&$K_{n,m}(f)(p)$\\
		\hline
		50& 50&1&0.15509756+4.9677361$\times 10^{-15}$i\\
		\hline
		100&100&1&0.15515498+2.968903571$\times 10^{-15}$i\\
		\hline
		200&200&1&0.155169346-1.47290502$\times 10^{-15}$i\\
		\hline
		300&300&1&0.1551720046+5.9769031$\times 10^{-15}$i\\
		\hline
		400&400&1&0.1551729368+3.46106193$\times 10^{-15}$i\\
		\hline
		500&500&1&0.1551733674+1.11732524$\times 10^{-16}$i\\
		\hline
	\end{tabular}
\end{center}

\begin{center}
	\begin{tabular}{ |p{2cm}|p{2cm}|p{1.5cm}|p{6 cm}| }
		\hline
		&&&\\
		$n$& $m$& p&$K_{n,m}(f)(p)$\\
		\hline
		50& 50&$\frac{\pi}{4}$&0.132503942-0.046439i\\
		\hline
		100&100&$\frac{\pi}{4}$&0.1325853-0.04646751i                        \\
		\hline
		200&200&$\frac{\pi}{4}$&0.13260566-0.0464746492i\\
		\hline
		300&300&$\frac{\pi}{4}$&0.13260943-0.0464759714i\\
		\hline
		400&400&$\frac{\pi}{4}$&0.132610751-0.04647643i\\
		\hline
		500&500&$\frac{\pi}{4}$&0.13261136-0.0464766i\\
		\hline
	\end{tabular}
\end{center}
\begin{center}
	\begin{tabular}{ |p{2cm}|p{2cm}|p{1.5cm}|p{6 cm}| }
		\hline
		&&&\\
		$n$& $m$& p&$K_{n,m}(f)(p)$\\
		\hline
		50& 50&1.5&0.095911405+0.0959114i\\
		\hline
		100&100&1.5&0.09583492+0.09583492i\\
		\hline
		200&200&1.5&0.09581573+0.09581573i\\
		\hline
		300&300&1.5&0.09581217+0.09581217i\\
		\hline
		400&400&1.5&0.09581093+0.09581093i\\
		\hline
		500&500&1.5&0.0958104+0.09581036i\\
		\hline
	\end{tabular}
\end{center}
Consider the extension $\{H_{n,\delta}\}_{n\in\mathbb{N},\delta>0}$ of $\{K_n\}_{n\in\mathbb{N}}$ to $L^1(\mathbb{R})$ (see Definition \ref{hndelta}). We obtain the following convergence result for $\{H_{n,\delta}\}_{n\in\mathbb{N},\delta>0}$.
\begin{theorem}\label{extn1}
	Suppose $f\in L^1(\mathbb{R})$. Let $\{\phi_\delta\}_{\delta>0}$ be an approximate identity such that $\hat{\phi_\delta}$ is compactly supported for each $\delta$. Then 
	\[
	\lim\limits_{\delta\to 0}\lim\limits_{n\to\infty}\|(H_{n,\delta}(f)-f)\chi_{[a,b]}\|_1= 0,
	\]
	for any compact interval $[a,b]$ of $\mathbb{R}$.
\end{theorem}
\begin{proof}
	Let $\{\phi_\delta\}_{\delta>0}$ be an approximate identity such that $\widehat{f*\phi_\delta}$ is compactly supported for each $\delta>0$.
	We have
	\begin{align*}
		H_{n,\delta}(f)(x)
		&=\frac{1}{2\pi}\sum_{l=-m_\delta}^{m_\delta}\int\limits_{l\pi}^{(l+1)\pi}\frac{1}{2}\sum\limits_{k=1}^n\widehat{f*\phi_\delta}(\theta_k^{(n)}+l\pi)(P_k^l(\theta+\frac{\pi}{2n})+P_k^l(\theta-\frac{\pi}{2n}))e^{ix\theta}\ d\theta\\
	\end{align*}
	for some $m_\delta\in\mathbb{N}$ which depends on $\phi_\delta$. Let $\|.\|_1$ be considered on some compact interval of $\mathbb{R}$. Consider,
	\begin{equation}\label{eq}
		\|H_{n,\delta}(f)-f\|_1\leq \|H_{n,\delta}(f)-f*\phi_\delta\|_1+\|f*\phi_\delta-f\|_1.
	\end{equation}
	As $n\rightarrow \infty$, the first term in the right hand side of \ref{eq} converges to $0$ by Remark \ref{rem}. As $\delta\to0$, the second term also converges to $0$ since $\{\phi_\delta\}_{\delta>0}$ is an approximate identity.
\end{proof}	
\begin{remark}
	In the definition of $H_{n,\delta}$, suppose that $\phi_\delta$ is continuous for each $\delta>0$. Let $x\in [0,\pi]$. We have,
	
	$|H_{n,\delta}(f\chi_{[a,b]})(x)|$
	\begin{equation*}
		\begin{split}
			&\leq \frac{1}{2\pi}\sum\limits_{l=-m_\delta}^{m_\delta}\frac{1}{2}\sum\limits_{k=1}^n|\widehat{f*\phi_\delta}(\theta_k^{(n)}+l\pi)|\int\limits_{l\pi}^{(l+1)\pi}|P_k^l(\theta+\frac{\pi}{2n})+P_k^l(\theta-\frac{\pi}{2n})|\ d\theta.\\
			&\leq \|\widehat{f\chi_{[0,\pi]}*\phi_\delta}\|_\infty \frac{1}{2\pi}\sum\limits_{l=-m_\delta}^{m_\delta}\frac{1}{2}\sum\limits_{k=1}^n\int\limits_{l\pi}^{(l+1)\pi}|P_k^l(\theta+\frac{\pi}{2n})+P_k^l(\theta-\frac{\pi}{2n})|\ d\theta\\
			&\leq \|f\chi_{[0,\pi]}\|_1\|\phi_\delta\|_\infty\frac{1}{2\pi}\sum\limits_{l=-m_\delta}^{m_\delta}c_1\pi \leq C_\delta\|f\|_1
		\end{split}
	\end{equation*}
	, where $C_\delta$ does not depend on $n$ or $f$ and $\|.\|_\infty$ denotes the sup-norm over $[-\pi,\pi]$. Hence, $H_{n,\delta}$ is a bounded linear operator on $L^1[0,\pi]$ for each $n\in \mathbb{N}$.
\end{remark}
\subsection{Rate of convergence}
In this section, we prove that on a dense subclass of $L^1(\mathbb{R})$, the extended operator $K_n$ has a uniform order of convergence.
\begin{theorem}\label{rate}
	Suppose $f\in L^1(\mathbb{R})$ such that $\widehat{f}$ is compactly supported. Then 
	\[
	\|(K_n(f)-f)\chi_{[a,b]}\|_1=\mathcal{O}(\xi_n)
	\]
	where $\xi_n=\frac{1}{2}\sup\limits_{\eta\in[0,\pi]}\sum\limits_{k=1}^n|\theta_k^{(n)}-\eta||P_k(\eta-\frac{\pi}{2n})+P_k(\eta+\frac{\pi}{2n})|=\mathcal{O}(\omega(\arccos,\frac{1}{n})) \ (n\to\infty)$, for any  compact interval $[a,b]$ of $\mathbb{R}$.
\end{theorem}
\begin{proof}
	Suppose that $\widehat{f}$ is compactly supported. We see that
	\[
	|K_n(f)(x)-f(x)|
	\leq\frac{1}{2\pi} \int\limits_{-\infty}^{\infty}|G_n(\widehat{f})(\theta)-\widehat{f}(\theta)|\ d\theta
	=\sum\limits_{l=-m}^{m}\int\limits_{l\pi}^{(l+1)\pi}|G_n(\widehat{f})(\theta)-\widehat{f}(\theta)|\ d\theta
	\]
	for some $m\in \mathbb{N}$.
	We also see that for $\theta\in[l\pi,(l+1)\pi]$ and  $\eta=\theta-l\pi$, $\eta\in [0,\pi]$. Now,
	\begin{multline*}
		\sum\limits_{k=1}^n|\widehat{f}(\theta_k^{(n)}+l\pi)-\widehat{f}(\theta)||P_k^l(\theta-\frac{\pi}{2n})+P_k(\theta+\frac{\pi}{2n})|=\\
		\sum\limits_{k=1}^n|\widehat{f}(\theta_k^{(n)}+l\pi)-\widehat{f}(\eta+l\pi)||P_k(\eta-\frac{\pi}{2n})+P_k(\eta+\frac{\pi}{2n})|.
	\end{multline*}
	By the properties of the modulus of continuity, for $\delta>0$,
	\[
	|\widehat{f}(\theta_k^{(n)}+l\pi)-\widehat{f}(\eta+l\pi)|\leq (1+\delta^{-1}|\theta_k^{(n)}-\theta|)\omega(\widehat{f},\delta)_l,
	\]
	for every $l\in \mathbb{Z}$. We denote $\omega(\widehat{f},\delta)_l$ to be the modulus of continuity of $\widehat{f}$ in the interval $[l\pi,(l+1)\pi]$.
	Suppose $\delta=\xi_n=\frac{1}{2}\sup\limits_{\eta\in[0,\pi]}\sum\limits_{k=1}^n|\theta_k^{(n)}-\eta||P_k(\eta-\frac{\pi}{2n})+P_k(\eta+\frac{\pi}{2n})|$. Finally we have,
	\begin{equation}\label{eqn1}
		\begin{split}
			|\int\limits_{-\infty}^{\infty}(G_n(\widehat{f})(\theta)-\widehat{f}(\theta))e^{ix\theta}\ d\theta|
			&\leq\sum\limits_{l=-m}^{m}\int\limits_{l\pi}^{(l+1)\pi}|G_n(\widehat{f})(\theta)-\widehat{f}(\theta)|\ d\theta\\
			&=\pi (c_1+1)\sum\limits_{l=-m}^{m}\omega(\widehat{f},\xi_n)_l.
		\end{split}
	\end{equation}
	Therefore, for $x\in\mathbb{R}$, $|\widecheck{G_n\widehat{(f)}}(x)-f(x)|\leq \pi (c_1+1)\sum\limits_{l=-m}^{m}\omega(\widehat{f},\xi_n)_l$.
	
	Let $x,y\in \mathbb{R}$. Then, 
	\begin{equation*}
		|\widehat{f}(x)-\widehat{f}(y)|=\sum\limits_{l=-m}^m\int\limits_{l\pi}^{(l+1)\pi}|\widehat{f}(\theta)||e^{-ix\theta}-e^{-iy\theta}|d\theta \leq C|x-y|,
	\end{equation*}
	since $|e^{-ix\theta}-e^{-iy\theta}|\leq \sqrt{2}|x-y||\theta|$ and $C=\sqrt{2}\sum\limits_{l=-m}^m\int\limits_{l\pi}^{(l+1)\pi}|\widehat{f}(\theta)||\theta|\ d\theta$. 
	For $|x-y|\leq \xi_n$, we have $|\widehat{f}(x)-\widehat{f}(y)|\leq M\xi_n$. Thus
	\begin{equation}\label{eqn2}
		\omega(\widehat{f},\xi_n)_l=\mathcal{O}(\xi_n)
	\end{equation}
	for every $l\in\mathbb{Z}$ where $-m\leq l\leq m$. Then by \ref{eqn1} and \ref{eqn2}, $\|(K_n(f)-f)\chi_{[a,b]}\|_1=\mathcal{O}(\xi_n)=\mathcal{O}(\omega(\arccos,\frac{1}{n}))\ (n\to\infty)$ (Theorem \ref{moc})
	for any compact interval $[a,b]$ of $\mathbb{R}$. 
\end{proof}

\begin{remark}\label{conv}
	Suppose $f\in L^1(\mathbb{R})$.	Note that $K_n(f*\phi_{\delta})=H_{n,\delta}(f)$. We also have that $\widehat{f*\phi_\delta}$ is compactly supported. Thus by the above result 
	\[
	\|H_{n,\delta}(f\chi_{[0,\pi]})-f\chi_{[0,\pi]}*\phi_\delta\|_1=\mathcal{O}(\xi_n),
	\]
	where $\|.\|_1$ is the $L_1$ norm on $[0,\pi]$. 
	We also have $\lim\limits_{\delta\to 0}\|H_{n,\delta}(f\chi_{[0,\pi]})-f\|_1=\mathcal{O}(\xi_n)$.
\end{remark}
We provide a numerical table that gives an idea about the convergence rate $\xi_n$. The above estimate indicates that, for any $f$ in the specified class, $K_n(f) \sim f $ with an error of $\mathcal{O}(\xi_n)$ in the $L^1$ sense.

\begin{center}
	\begin{tabular}{ |p{2cm}|p{4cm}|  }
		\hline
		&\\
		n& $\xi_n$\\
		\hline
		
		100&0.04436868245\\
		\hline
		200&0.01970873322\\
		\hline
		300&0.01513454629\\
		\hline
		400&0.01448400327\\
		\hline
		500&0.011950746336\\
		\hline
		600&0.010190080703\\
		\hline
		700&0.009029628692\\
		\hline
		800&0.008142319697\\
		\hline
		900&0.007361966592\\
		\hline
		1000&0.003879745159\\
		\hline
	\end{tabular}
\end{center}
\subsection*{Concluding Remarks and Future Problems}
In this article, we proved an operator analogue of Wulbert's theorem (Theorem \ref{wul1}), which is a non-positive version of Korovkin theorem. We also derived an analogous result on the space $L^1[0,1]$, under additional assumptions. A more general formulation in this setting remains an open direction for further research. We illustrated these results on concrete examples. Extending these results to the general case on $L^p[0,1]$, $1<p<+\infty$ spaces is a research scope.

\subsection*{Acknowledgements}
V. B. Kiran Kumar is supported by the KSYSA-Research Grant by KSCSTE, Kerala. P. C. Vinaya is thankful the University Grants Commission (UGC), Government of India for financial support.


\end{document}